# HOPF POTENTIALS FOR THE SCHRÖDINGER OPERATOR

LUIGI ORSINA AND AUGUSTO C. PONCE

ABSTRACT. We establish the Hopf boundary point lemma for the Schrödinger operator $-\Delta + V$ involving potentials $V$ that merely belong to the space $L^1_{\text{loc}}(\Omega)$. More precisely, we prove that among all nonnegative supersolutions $u$ of $-\Delta + V$ which vanish on the boundary $\partial \Omega$ and are such that $Vu \in L^1(\Omega)$, if there exists *one* supersolution that satisfies $\partial u/\partial n < 0$ almost everywhere on $\partial \Omega$ with respect to the outward unit vector $n$, then such a property holds for *every* nontrivial supersolution in the same class. We rely on the existence of nontrivial solutions of the nonhomogeneous Dirichlet problem with boundary datum in $L^\infty(\partial \Omega)$.

## 1. INTRODUCTION AND MAIN RESULTS

Let $\Omega \subset \mathbb{R}^N$ be a smooth bounded connected open set. The Hopf boundary point lemma for elliptic PDEs asserts that if $u \in C^2(\Omega) \cap C^1(\overline{\Omega})$ satisfies the Dirichlet problem

(1.1) $$\begin{cases} -\Delta u = \mu & \text{in } \Omega, \\ u = 0 & \text{on } \partial\Omega, \end{cases}$$

and if $\mu \geq 0$ in $\Omega$, then the normal derivative of $u$ with respect to the outward unit vector $n$ satisfies

$$\frac{\partial u}{\partial n} < 0 \quad \text{on } \partial\Omega;$$

see [25, Proposition A.4.1; 26, Section 6.4.2; 29, Lemma 3.4]. The classical weak maximum principle states that $u \geq 0$ in $\Omega$, so the information that $\partial u/\partial n \leq 0$ merely follows from the minimality of $u$ on $\partial\Omega$. The main issue involving the Hopf lemma is that if $\partial u(a)/\partial n = 0$ for some $a \in \partial\Omega$, then $u \equiv 0$ in $\Omega$. A drawback of this formulation lies in the $C^1$ regularity of $u$ that is required near the boundary.

The Hopf lemma can be also stated quantitatively, based on the Morel–Oswald maximum principle as follows:

(1.2) $$u(x) \geq c \left( \int_\Omega \mu \, d_{\partial\Omega} \right) d_{\partial\Omega}(x),$$

---

2010 *Mathematics Subject Classification.* Primary: 35B05, 35B50; Secondary: 31B15, 31B35.

where $c > 0$ and $d_{\partial\Omega} : \overline{\Omega} \to \mathbb{R}$ denotes the distance to the boundary; see [12, Lemma 3.2; 25, Proposition A.4.2]. This inequality is equivalent to the pointwise estimate for the Green's function [48]:

$$G(x, y) \geq c' \, d_{\partial\Omega}(x) d_{\partial\Omega}(y).$$

By smoothness of the solution $\zeta$ of the Dirichlet problem (1.1) with constant density $\mu \equiv 1$, we have $d_{\partial\Omega} \geq a\zeta$ for some constant $a > 0$. Inserting this estimate in (1.2) and integrating by parts, we can thus rewrite inequality (1.2) as

$$(1.3) \qquad u(x) \geq c'' \left( \int_\Omega u \right) d_{\partial\Omega}(x).$$

Since (1.2) and (1.3) do not explicitly involve the normal derivative, by an approximation argument both inequalities also apply to the solutions of the Dirichlet problem involving rougher data $\mu$, for instance in the class of $L^1$ functions or finite measures.

Another effective approach to avoid smoothness of the solution near the boundary consists in associating a notion of distributional normal derivative when $\mu$ is merely a finite measure in $\Omega$, and the equation is satisfied in the sense of distributions:

$$(1.4) \qquad -\int_\Omega u \, \Delta\varphi = \int_\Omega \varphi \, \mathrm{d}\mu,$$

for every test function $\varphi \in C_c^\infty(\Omega)$ with compact support in $\Omega$. The zero boundary datum can be encoded by assuming that $u \in W_0^{1,1}(\Omega)$. An equivalent strategy to give a meaning to the Dirichlet problem consists in using the integral formulation (1.4) with the larger class $C_0^\infty(\overline{\Omega})$ of test functions which are smooth in $\overline{\Omega}$ and vanish on $\partial\Omega$; see [34; 42, Proposition 6.3].

It has been observed by Brezis and the second author in [15, Theorem 1.2] that if $u \in W_0^{1,1}(\Omega)$ satisfies (1.4) for some finite measure $\mu$ in $\Omega$, then there exists a unique function $F \in L^1(\partial\Omega)$ such that

$$(1.5) \qquad \int_\Omega \nabla u \cdot \nabla \psi = \int_\Omega \psi \, \mathrm{d}\mu + \int_{\partial\Omega} \psi F \, \mathrm{d}\sigma,$$

for every test function $\psi \in C^\infty(\overline{\Omega})$, where $\sigma = \mathcal{H}^{N-1}\lfloor_{\partial\Omega}$ denotes the surface measure of $\partial\Omega$. The distributional normal derivative of $u$ is then defined as

$$\frac{\partial u}{\partial n} := F.$$

When $u$ is smooth on $\overline{\Omega}$, the notions of classical and distributional normal derivatives of $u$ coincide. The definition of a distributional pairing $\langle \partial u/\partial n, \psi \rangle$ under various assumptions on $u$ and $\nabla u$ has been investigated by several authors; see e.g. [5, 32, 33]. The main point here is its realization as a legitimate function in $L^1(\partial\Omega)$, like in [6]. We refer the reader to [3, 15] for additional properties of the distributional normal derivative, including the setting of nonhomogeneous Dirichlet problems.

Identity (1.5) implies in particular that the function $u$, extended by zero to $\mathbb{R}^N$, is such that $\Delta u$ is a finite measure in $\mathbb{R}^N$. The distributional normal derivative satisfies a comparison estimate (see Proposition 2.1 below) which combined with (1.3) provides one with the uniform bound

$$\text{(1.6)} \qquad \frac{\partial u}{\partial n}(x) \leq -c'' \int_\Omega u.$$

One of the motivations of our work comes from the seminal paper of Kato [31] on the Schrödinger operator $-\Delta + V$ involving potentials $V$ which are merely $L^1_{\text{loc}}(\Omega)$. The Hopf lemma above has an affirmative counterpart for potentials $V \in L^\infty(\Omega)$, but we are interested in situations where $V$ need not be summable near the boundary.

Many classical properties that hold for the Laplacian need no longer be true for $-\Delta + V$ due to some possible singular behavior of $V$. In this regard, two instructive examples are provided by the smooth functions $u_i : \overline{B}_1 \to \mathbb{R}$ defined by

$$u_1(x) = (1 - |x|^2)^2 \quad \text{and} \quad u_2(x) = (1 - |x|^2)|x - a|,$$

where $B_1 := B_1(0)$ denotes the unit ball in $\mathbb{R}^N$ centered at the origin and $a \in \partial B_1$ is any given point on the boundary. In the first case, we have that $\partial u_1 / \partial n \equiv 0$ on $\partial B_1$ and the Schrödinger equation

$$-\Delta u_1 + V u_1 = 0 \quad \text{in } B_1$$

is satisfied in terms of a potential $V$ that behaves like $1/d_{\partial B_1}^2$ near the boundary; in particular, $V \notin L^1(B_1; d_{\partial B_1} \, \mathrm{d}x)$. In the second case, we have that $\partial u_2 / \partial n < 0$ except at $a$ and the Schrödinger equation is now satisfied for another potential $V$ such that $V \in L^1(B_1; d_{\partial B_1} \, \mathrm{d}x)$ in dimension $N \geq 2$.

We thus have the appearance of an exceptional set where the Hopf lemma fails, and it is our goal in this paper to understand how big such an exceptional set can be. A more refined example which we develop in Section 7 below shows that *every* compact subset $K \subset \partial \Omega$ with zero surface measure is an exceptional set for some suitable potential $V \in L^1(\Omega; d_{\partial \Omega} \, \mathrm{d}x)$. It follows from our Theorem 1 below that there can be essentially no other exceptional sets in this case.

The class of functions we consider consists of supersolutions $u \in W^{1,1}_0(\Omega)$ of the Schrödinger operator $-\Delta + V$ in the sense of distributions. More precisely, we assume that $Vu \in L^1_{\text{loc}}(\Omega)$ and

$$\int_\Omega u \left( -\Delta \varphi + V \varphi \right) \geq 0,$$

for every nonnegative function $\varphi \in C^\infty_c(\Omega)$. By a classical property in the Theory of Distributions, we have in this case that $-\Delta u + Vu$ is a locally finite measure in $\Omega$. However, the measure $\Delta u$ need not be finite in $\Omega$ and so the distributional normal derivative may not be well-defined in $L^1(\partial \Omega)$ in the sense of [15]. For this reason, we define in this paper the normal derivative for functions $u$ that are merely in $W_0^{1,1}(\Omega)$: by $\partial u/\partial n$ we mean the essential infimum of the set

$$\Big\{\frac{\partial w}{\partial n} \in L^1(\partial\Omega) : w \in \mathcal{G}_u\Big\},$$

where $\mathcal{G}_u$ denotes the class of functions $w \in W_0^{1,1}(\Omega)$ such that $\Delta w$ is a finite measure in $\Omega$ and $w \le u$ almost everywhere in $\Omega$; see Section 2 below. In particular, if $u$ is nonnegative, then the normal derivative $\partial u/\partial n$ is a Borel function with values in $[-\infty, 0]$, and if we happen to know that $\Delta u$ is a finite measure in $\Omega$, then $u \in \mathcal{G}_u$ and $\partial u/\partial n$ coincides with the distributional normal derivative.

To understand the mechanism that is hidden behind the examples above concerning the failure of the Hopf lemma, we introduce the concept of Hopf potential as follows:

**Definition 1.1.** *We say that $V \in L^1_{\text{loc}}(\Omega)$ is a* Hopf potential *whenever there exists a nonnegative function $\zeta_0 \in W_0^{1,1}(\Omega)$ such that*

($H_1$) $V\zeta_0 \in L^1(\Omega)$,
($H_2$) $\partial\zeta_0/\partial n < 0$ *almost everywhere on $\partial\Omega$*.

As a trivial consequence of this definition, for every Hopf potential $V$ and every $\alpha \in \mathbb{R}$, the function $\alpha V$ is also a Hopf potential. We show in Section 2 that the class of Hopf potentials is actually a vector subspace of $L^1_{\text{loc}}(\Omega)$. Since the solution $\zeta$ of the Dirichlet problem (1.1) with constant density $\mu \equiv 1$ behaves as $d_{\partial\Omega}$ near the boundary by the classical Hopf lemma, we have that $V\zeta \in L^1(\Omega)$ if and only if $V \in L^1(\Omega; d_{\partial\Omega}\,\mathrm{d}x)$. Therefore, *every $V \in L^1(\Omega; d_{\partial\Omega}\,\mathrm{d}x)$ is a Hopf potential.*

We establish the following qualitative counterpart of estimate (1.6) for $-\Delta + V$ when $V$ is a Hopf potential:

**Theorem 1.** *Let $V \in L^1_{\text{loc}}(\Omega)$ be a Hopf potential and let $u \in W_0^{1,1}(\Omega)$ be a nonnegative supersolution of the Schrödinger operator $-\Delta + V$. If $Vu \in L^1(\Omega)$ and $\int_\Omega u > 0$, then*

$$\frac{\partial u}{\partial n} < 0 \quad \text{almost everywhere on } \partial\Omega.$$

Theorem 1 above contains as a particular case a Hopf lemma by Bertsch, Smarrazzo and Tesei [10, Proposition 3.4] which implies the main result in their paper (Theorem 2.1) concerning a characterization of the strong maximum principle in dimension $N = 1$; see also [9, Lemma 3.6]. To tackle the Hopf lemma in any dimension $N \ge 1$, we rely on a different strategy based on a careful combination of fine properties from Measure theory and Elliptic PDEs.

One may also consider a localized counterpart of the concept of Hopf potentials, where property ($H_2$) need not be satisfied by $\zeta_0$ on the entire boundary, but only on a subset of it. In fact, we deduce Theorem 1 from a more general result which is valid for potentials $V$ that merely belong to $L^1_{\mathrm{loc}}(\Omega)$:

**Theorem 2.** *Let $V \in L^1_{\mathrm{loc}}(\Omega)$ and let $u_i \in W^{1,1}_0(\Omega)$, with $i \in \{1,2\}$, be two nonnegative supersolutions of the Schrödinger operator $-\Delta + V$. If $Vu_i \in L^1(\Omega)$ and $\int_\Omega u_i > 0$, then for almost every $x \in \partial\Omega$ we have*

$$\frac{\partial u_1}{\partial n}(x) < 0 \quad \textit{if and only if} \quad \frac{\partial u_2}{\partial n}(x) < 0.$$

This theorem yields the remarkable property that once there exists *one* supersolution for $-\Delta + V$ satisfying the conclusion of the classical Hopf lemma on a subset $A \subset \partial\Omega$, then *every* supersolution also satisfies the Hopf lemma on $A$ except for a negligible subset of $\partial\Omega$. Such a conclusion bears some striking analogy with the (straightforward) generalized weak maximum principle for linear elliptic operators of second order [44, Chapter 2, Theorem 10]: for the Schrödinger operator $-\Delta + V$ with a possibly signed potential $V$, the existence of one positive supersolution implies that every nonzero supersolution which vanishes on the boundary must be positive.

The existence of a positive supersolution is also equivalent to the positivity of the energy functional

$$\varphi \in C^\infty_c(\Omega) \longmapsto \int_\Omega (|\nabla\varphi|^2 + V\varphi^2),$$

which is at the heart of the Agmon–Allegretto–Piepenbrinck positivity principle; see e.g. [25, Theorem A.6.1] and also [21, 39] for more detailed information and further perspectives. Although one typically assumes that $V \in L^p_{\mathrm{loc}}(\Omega)$ with $p > N/2$, the validity of the strong maximum principle when $V \in L^p_{\mathrm{loc}}(\Omega)$ with $p \geq 1$ (see [1, 38]) support for an extension of the Agmon–Allegretto–Piepenbrinck principle for potentials $V$ that merely belong to $L^1_{\mathrm{loc}}(\Omega)$ based on the tools we develop to prove Theorems 1 and 2 above; see [8, 41] for Hardy potentials and [40] for potentials in Morrey spaces.

A key ingredient of our analysis relies on Proposition 5.1 below which establishes the equivalence between the validity of the Hopf lemma for the Schrödinger operator $-\Delta + V$ and the existence of nontrivial solutions of the nonhomogeneous Dirichlet problem

(1.7) $$\begin{cases} -\Delta w + Vw = 0 & \text{in } \Omega, \\ \phantom{-\Delta w + V}w = g & \text{on } \partial\Omega, \end{cases}$$

with nonnegative potentials $V \in L^1_{\mathrm{loc}}(\Omega)$, for any datum $g \in L^\infty(\partial\Omega)$. The meaning of a solution of (1.7) is a delicate issue due to the possible singular behavior of $V$ near the boundary. Our approach is based on the use of nonsmooth test functions that satisfy a Dirichlet problem involving interior measure data in the spirit of Stampacchia's definition of weak solutions via duality [45]. That problem (1.7) has a solution in this sense for every $g \in L^\infty(\partial\Omega)$ can be handled using an approximation procedure starting from variational solutions; see Section 3.

Our strategy to tackle (1.7) differs from the recent work of Véron and Yarur [47] that investigates problem (1.7) with finite boundary measure data and nonnegative potentials $V \in L^\infty_{\text{loc}}(\Omega)$. They rely on the definition of a solution using test functions like $C_0^\infty(\overline{\Omega})$, which do not take into account the singular behavior of the potential $V$, and on the Poisson representation of the solution in terms of the Poisson kernel associated to $-\Delta + V$.

Due to the singular behavior of $V$, it may happen in our case that $w \equiv 0$ is the (unique) solution of (1.7) even if $g \neq 0$. An example of such a counterintuitive phenomenon is given by any potential $V \sim 1/d_{\partial\Omega}^2$, for which the Hopf lemma fails completely. In this case, the equation

$$-\Delta w + Vw = 0 \quad \text{in } \Omega$$

can have nontrivial solutions but they satisfy some normalized boundary trace that has been investigated by Marcus and Nguyen [35].

Another strategy that has been pursued by Ancona [2] is based on the existence of the Martin kernel $K_a^V$ for $a \in \partial\Omega$ under the assumption that $V$ is a potential in $L^\infty_{\text{loc}}(\Omega)$ that satisfies

$$(1.8) \qquad 0 \leq V \lesssim \frac{1}{d_{\partial\Omega}^2}.$$

For instance, in the setting of positive solutions of the semilinear equation

$$-\Delta u + u^q = 0 \quad \text{in } \Omega$$

with exponent $q > 1$, the potential $V = u^{q-1}$ satisfies (1.8) by the Keller–Osserman estimate; see [36, Chapter 4]. In general, the study of fine regular points of the Schrödinger operator $-\Delta + V$ through Martin kernels gives another approach to the existence of solutions of (1.7). In this regard, Ancona [4, 47] proved that $a \in \partial\Omega$ is a fine regular point for $-\Delta + V$ if and only if

$$(1.9) \qquad \int_\Omega \frac{d_{\partial\Omega}^2(x)}{|x-a|^N} V(x)\,\mathrm{d}x < +\infty.$$

When, in addition to (1.8), $V$ belongs to $L^1(\Omega; d_{\partial\Omega}\,\mathrm{d}x)$, integration of the left-hand side of (1.9) over $\partial\Omega$ with respect to $a$ and Fubini's theorem imply that almost every $a \in \partial\Omega$ satisfies (1.9). This agrees with our conclusion concerning the existence of nontrivial solutions for (1.7) since we know in this case that $V$ is a Hopf potential. It is unclear however how one can avoid assumption (1.8) in this setting: Ancona's argument strongly relies on the Harnack principle, which is not true when one merely has $V \in L^1_{\text{loc}}(\Omega)$.

Observe that from the physical point of view the infinite-potential well $1/d_{\partial\Omega}^2$ is so strong that it confines particles inside $\Omega$, which mathematically means that supersolutions must have a vanishing normal derivative on $\partial\Omega$; see [22–24] and also Example 8.2 below. Although such a conclusion can be successfully deduced from Theorem 2 by looking explicitly for one supersolution such that $\partial u/\partial n \equiv 0$ on $\partial\Omega$, we give a direct proof of this fact by a simple measure-theoretic argument that does not rely on the PDE; see Proposition 2.7.

The paper is organized as follows. In Section 2 we extend the concept of normal derivative to any function in $W_0^{1,1}(\Omega)$, even if $\Delta u$ is not a finite measure in $\Omega$. In Section 3, we prove the existence of solutions of the nonhomogeneous Dirichlet problem with $L^\infty$ data; the meaning of solution is given by means of duality with solutions of the Dirichlet problem with measure data. In Section 4, we prove the existence of nonnegative solutions of the nonhomogeneous problem when the boundary datum is nonnegative but the inner datum is nonpositive. We then explain how this property implies Theorem 1 in the case of smooth supersolutions. In Section 5 we explain the connection between the Hopf lemma and the existence of nontrivial solutions of (1.7). Theorems 1 and 2 are then proved in Section 6. We then show in Section 7 that every negligible compact subset of $\partial\Omega$ is the zero-set $\{\partial u/\partial n = 0\}$ for some smooth positive solution of the Schrödinger equation $-\Delta u + Vu = 0$ such that $V \in L^1(\Omega; d_{\partial\Omega}\,\mathrm{d}x)$. In Section 8 we explain why Theorems 1 and 2 cannot be true for potentials $V: \Omega \to [0, +\infty]$ that are merely Borel functions.

## 2. Normal derivative as a Borel function

The notion of distributional normal derivative from [15] applies to any function $u \in W_0^{1,1}(\Omega)$ such that $\Delta u$ is a finite measure in $\Omega$. In this case, the normal derivative $\partial u/\partial n$ is an element in $L^1(\partial\Omega)$ such that

$$\int_\Omega \nabla u \cdot \nabla \psi = -\int_\Omega \psi\,\Delta u + \int_{\partial\Omega} \psi\,\frac{\partial u}{\partial n}\,\mathrm{d}\sigma \quad \text{for every } \psi \in C^\infty(\overline{\Omega}).$$

In this work, we deal with functions $u \in W_0^{1,1}(\Omega)$ such that the distribution $\Delta u$ need not be a finite measure. The strategy we adopt to define a Borel normal derivative is motivated by the following comparison principle which can be deduced from Kato's inequality (see [42, Lemma 12.15]):

**Proposition 2.1.** *Let $v \in W_0^{1,1}(\Omega)$ be such that $\Delta v$ is a finite measure in $\Omega$. If $v \geq 0$ almost everywhere in $\Omega$, then $\partial v/\partial n \leq 0$ almost everywhere on $\partial\Omega$ with respect to the surface measure.*

Now, to our definition of normal derivative as a Borel function, we begin with any $u \in W_0^{1,1}(\Omega)$. By the essential infimum of the set

$$\mathcal{N}_u := \Big\{ \frac{\partial w}{\partial n} \in L^1(\partial\Omega) : w \in \mathcal{G}_u \Big\},$$

where

$$\mathcal{G}_u := \Big\{ w \in W_0^{1,1}(\Omega) : \Delta w \text{ is a finite measure and } w \le u \text{ a.e. in } \Omega \Big\},$$

we mean a Borel function $F : \partial\Omega \to [-\infty, \infty]$ such that

($i$) $F \le \partial w/\partial n$ almost everywhere on $\partial\Omega$, for every $w \in \mathcal{G}_u$,
($ii$) if $\widetilde{F} : \partial\Omega \to [-\infty, \infty]$ is another Borel function that satisfies ($i$), then $\widetilde{F} \le F$ almost everywhere on $\partial\Omega$.

We then define the normal derivative of $u$ as

$$\frac{\partial u}{\partial n} := F.$$

**Proposition 2.2.** *Such a normal derivative $\partial u/\partial n$ exists for every $u \in W_0^{1,1}(\Omega)$.*

*Proof.* By the separability of $L^1(\partial\Omega)$, we can extract a countable subset $A$ of $\mathcal{G}_u$ such that $\{\partial v/\partial n : v \in A\}$ is dense in $\mathcal{N}_u$. We claim that the Borel measurable function $F : \partial\Omega \to [-\infty, \infty]$ defined by

$$F(x) := \inf_{v \in A} \frac{\partial v}{\partial n}(x)$$

satisfies Properties ($i$) and ($ii$) above. Indeed, given $w \in \mathcal{G}_u$, take a sequence $(v_k)_{k \in \mathbb{N}}$ in $A$ such that $(\partial v_k/\partial n)_{k \in \mathbb{N}}$ converges to $\partial w/\partial n$ in $L^1(\partial\Omega)$. Passing to a subsequence if necessary, we may assume that the convergence holds almost everywhere on $\partial\Omega$. Since $\partial v_k/\partial n \ge F$ on $\partial\Omega$, we deduce that $\partial w/\partial n \ge F$ almost everywhere on $\partial\Omega$. Hence, $F$ satisfies Property ($i$). We now let $\widetilde{F}$ be another function that satisfies Property ($i$), and for each $v \in A$ denote by $E_v \subset \partial\Omega$ a set of surface measure zero such that $\widetilde{F}(x) \le \partial v(x)/\partial n$ for every $x \in \partial\Omega \setminus E_v$. Since $A$ is countable, the set $E = \bigcup_{v \in A} E_v$ also has surface measure zero and

$$\widetilde{F}(x) \le \frac{\partial v}{\partial n}(x) \quad \text{for every } x \in \partial\Omega \setminus E,$$

for every $v \in A$. Taking the infimum of the right-hand side over $v$ we deduce that $\widetilde{F} \le F$ on $\partial\Omega \setminus E$, which gives Property ($ii$). $\square$

In the pointwise approximation of a Borel normal derivative, one can restrict the attention to the study of monotone sequences in $\mathcal{G}_u$ and $\mathcal{N}_u$:

**Proposition 2.3.** *For every $u \in W_0^{1,1}(\Omega)$, there exists a nondecreasing sequence $(w_k)_{k \in \mathbb{N}}$ in $\mathcal{G}_u$ such that $(\partial w_k/\partial n)_{k \in \mathbb{N}}$ is a non-increasing sequence in $\mathcal{N}_u$ that converges almost everywhere to $\partial u/\partial n$ on $\partial\Omega$.*

In order to prove Proposition 2.3 we rely on Kato's inequality up to the boundary which implies that if $\zeta \in W_0^{1,1}(\Omega)$ and $\Delta \zeta$ is a finite measure in $\Omega$, then $\Delta[(\zeta - a)^+]$ is also a finite measure in $\Omega$ for every $a \in \mathbb{R}$ and

$$\|\Delta[(\zeta - a)^+]\|_{\mathcal{M}(\Omega)} \leq 2\|\Delta\zeta\|_{\mathcal{M}(\Omega)}; \tag{2.1}$$

see [15, Theorem 1.1; 42, Proposition 7.7]. Here, $\mathcal{M}(\Omega)$ denotes the vector space of finite Borel measures $\nu$ in $\Omega$ equipped with the norm

$$\|\nu\|_{\mathcal{M}(\Omega)} = |\nu|(\Omega),$$

that makes $\mathcal{M}(\Omega)$ a Banach space. The normal derivative $\partial(\zeta - a)^+/\partial n$ is then well-defined in the distributional sense as an element in $L^1(\partial\Omega)$. Ancona subsequently proved using tools from Potential theory that

$$\frac{\partial(\zeta - a)^+}{\partial n} = \begin{cases} \partial\zeta/\partial n & \text{if } a < 0, \\ \min\{\partial\zeta/\partial n, 0\} & \text{if } a = 0, \\ 0 & \text{if } a > 0; \end{cases} \tag{2.2}$$

see [3, Remark 6.2]. These properties can be illustrated by the following lemma:

**Lemma 2.4.** *If $\zeta_i \in W_0^{1,1}(\Omega)$, with $i \in \{1, 2\}$, are such that $\Delta\zeta_i$ are finite measures in $\Omega$, then the function $\zeta = \max\{\zeta_1, \zeta_2\}$ belongs to $W_0^{1,1}(\Omega)$, is such that $\Delta\zeta$ is a finite measure in $\Omega$, and*

$$\frac{\partial \zeta}{\partial n} = \min\left\{\frac{\partial \zeta_1}{\partial n}, \frac{\partial \zeta_2}{\partial n}\right\} \quad \text{almost everywhere on } \partial\Omega.$$

*Proof of Lemma 2.4.* Observe that

$$\zeta = \zeta_2 + (\zeta_1 - \zeta_2)^+.$$

Thus, $\zeta \in W_0^{1,1}(\Omega)$. By Kato's inequality up to the boundary (2.1) applied to the function $\zeta_1 - \zeta_2$ and $a = 0$, we deduce that the measure $\Delta\big[(\zeta_1 - \zeta_2)^+\big]$ is finite in $\Omega$, whence so is the measure $\Delta\zeta$. By (2.2) we have

$$\frac{\partial}{\partial n}(\zeta_1 - \zeta_2)^+ = \min\left\{\frac{\partial}{\partial n}(\zeta_1 - \zeta_2), 0\right\} = -\frac{\partial \zeta_2}{\partial n} + \min\left\{\frac{\partial \zeta_1}{\partial n}, \frac{\partial \zeta_2}{\partial n}\right\},$$

and the conclusion follows. $\square$

*Proof of Proposition 2.3.* Let $(v_k)_{k \in \mathbb{N}}$ be a sequence in $\mathcal{G}_u$ such that $(\partial v_k/\partial n)_{k \in \mathbb{N}}$ is dense in $\mathcal{N}_u$. As in the proof of Proposition 2.2, we have

$$\frac{\partial u}{\partial n} = \inf_{j \in \mathbb{N}} \frac{\partial v_j}{\partial n} \quad \text{almost everywhere on } \partial\Omega.$$

Define by induction the nondecreasing sequence $(w_k)_{k \in \mathbb{N}}$ as $w_0 := v_0$ and, for $k \in \mathbb{N}_*$,

$$w_k := \max\{w_{k-1}, v_k\}.$$

By Lemma 2.4 we have $w_k \in \mathcal{G}_u$ for every $k \in \mathbb{N}$. In particular, $\partial u/\partial n \leq \partial w_k/\partial n$ almost everywhere on $\partial\Omega$. By comparison of normal derivatives, the sequence $(\partial w_k/\partial n)_{k \in \mathbb{N}}$ is monotone and non-increasing, hence

$$\frac{\partial u}{\partial n} \leq \lim_{k \to \infty} \frac{\partial w_k}{\partial n}$$

and also

$$\lim_{k \to \infty} \frac{\partial w_k}{\partial n} \leq \frac{\partial w_j}{\partial n} \leq \frac{\partial v_j}{\partial n}$$

almost everywhere on $\partial\Omega$, for every $j \in \mathbb{N}$. Taking the infimum in the right-hand side over $j$, we deduce that

$$\lim_{k \to \infty} \frac{\partial w_k}{\partial n} = \frac{\partial u}{\partial n} \quad \text{almost everywhere on } \partial\Omega. \qquad \square$$

As a consequence of Proposition 2.3, we observe that for the sake of investigating the set where the normal derivative of a function $u$ is negative, one does not need to rely on the entire family $\mathcal{G}_u$ nor even on a countable subset of it, but on a single suitably-chosen element:

**Proposition 2.5.** *For every nonnegative function $u \in W_0^{1,1}(\Omega)$, there exists a nonnegative function $v \in \mathcal{G}_u$ such that*

$$\frac{\partial v}{\partial n} < 0 \quad \text{almost everywhere on } \Big\{\frac{\partial u}{\partial n} < 0\Big\}.$$

*Proof.* Let $(w_k)_{k \in \mathbb{N}}$ be a nondecreasing sequence in $\mathcal{G}_u$ satisfying the conclusion of Proposition 2.3. Replacing each $w_k$ by its positive part if necessary, we may assume by Lemma 2.4 that each function $w_k$ is nonnegative in $\Omega$ and in particular $\partial w_0/\partial n$ is nonpositive on $\partial\Omega$. Hence, in the sum

$$\frac{\partial w_0}{\partial n} + \sum_{j=1}^{\infty} \Big(\frac{\partial w_j}{\partial n} - \frac{\partial w_{j-1}}{\partial n}\Big) = \frac{\partial u}{\partial n},$$

we have $\partial w_0/\partial n \leq 0$ and $\partial w_j/\partial n - \partial w_{j-1}/\partial n \leq 0$ almost everywhere on $\partial\Omega$, for every $j \in \mathbb{N}_*$. In addition, for almost every $x \in \{\partial u/\partial n < 0\}$ one of these terms, possibly depending on $x$, must be negative. The conclusion is thus satisfied with

$$v := w_0 + \sum_{j=1}^{\infty} \epsilon_j (w_j - w_{j-1}),$$

where $(\epsilon_j)_{j \in \mathbb{N}}$ is a sequence in $(0, 1]$ such that

$$\sum_{j=1}^{\infty} \epsilon_j \big(\|\nabla(w_j - w_{j-1})\|_{L^1(\Omega)} + \|\Delta(w_j - w_{j-1})\|_{\mathcal{M}(\Omega)}\big) < +\infty.$$

Indeed, we have $v \leq u$ in $\Omega$ and such a choice of sequence $(\epsilon_j)_{j \in \mathbb{N}}$ ensures that $v$ belongs to $W_0^{1,1}(\Omega)$ and $\Delta v$ is a finite measure in $\Omega$ by completeness of $W_0^{1,1}(\Omega)$ and $\mathcal{M}(\Omega)$. $\qquad \square$

**Corollary 2.6.** *The class of Hopf potentials is a vector subspace of $L_{\text{loc}}^1(\Omega)$.*

*Proof.* Let $V_i \in L^1_{\text{loc}}(\Omega)$, with $i \in \{1, 2\}$, be two Hopf potentials, and denote by $\zeta_i \in W^{1,1}_0(\Omega)$ a nonnegative function that satisfies Properties $(H_1)$ and $(H_2)$ with respect to $V_i$. We now verify that $\alpha_1 V_1 + \alpha_2 V_2$ is a Hopf potential for every $\alpha_1, \alpha_2 \in \mathbb{R}$ by using the function $\zeta := \min\{\zeta_1, \zeta_2\}$. Observe that

$$\left|(\alpha_1 V_1 + \alpha_2 V_2)\zeta\right| \leq |\alpha_1 V_1|\zeta_1 + |\alpha_2 V_2|\zeta_2 \in L^1(\Omega).$$

By Proposition 2.5 there exists a nonnegative function $v_i \in \mathcal{G}_{\zeta_i}$ such that $\partial v_i/\partial n < 0$ almost everywhere on $\partial\Omega$. Since $\zeta \geq \min\{v_1, v_2\}$, the counterpart of Lemma 2.4 for the minimum of two functions gives in this case $\min\{v_1, v_2\} \in \mathcal{G}_\zeta$ and

$$\frac{\partial \zeta}{\partial n} \leq \frac{\partial}{\partial n}\min\{v_1, v_2\} = \max\left\{\frac{\partial v_1}{\partial n}, \frac{\partial v_2}{\partial n}\right\}.$$

Therefore, $\partial \zeta/\partial n < 0$ almost everywhere on $\partial\Omega$. $\square$

In the Schrödinger operator $-\Delta + V$, the behavior of the potential $V$ as $1/d^2_{\partial\Omega}$ near the boundary is critical for the Hopf lemma:

**Proposition 2.7.** *If $u \in W^{1,1}_0(\Omega)$ is a nonnegative function such that*

$$\tag{2.3} \int_\Omega \frac{u}{d^2_{\partial\Omega}} < \infty,$$

*then $\partial u/\partial n = 0$ almost everywhere on $\partial\Omega$.*

*Proof.* By the comparison principle for normal derivatives (Proposition 2.1), it suffices to verify that for every nonnegative function $v \in \mathcal{G}_u$ we have $\partial v/\partial n = 0$ almost everywhere on $\partial\Omega$. To this end, denote by $\mu$ the measure in $\mathbb{R}^N$ such that $\mu = \Delta v$ in $\Omega$ and $\mu \equiv 0$ on the Borel subsets of $\mathbb{R}^N \setminus \Omega$; we also extend $v$ by zero to $\mathbb{R}^N \setminus \Omega$. For every $\psi \in C^\infty(\mathbb{R}^N)$, we then have

$$\int_{\partial\Omega} \psi \frac{\partial v}{\partial n}\,\mathrm{d}\sigma = \int_\Omega \psi\,\Delta v + \int_\Omega \nabla v \cdot \nabla \psi = \int_{\mathbb{R}^N} \psi\,\mathrm{d}\mu - \int_{\mathbb{R}^N} v\,\Delta\psi.$$

Given $x \in \partial\Omega$ and $r > 0$, we apply this identity with $\psi(y) = \varphi((y - x)/r)$, where $\varphi \in C^\infty_c(\mathbb{R}^N)$ is such that $\varphi \equiv 1$ in $B_1$, $\varphi \equiv 0$ in $\mathbb{R}^N \setminus B_2$, and $0 \leq \varphi \leq 1$ in $\mathbb{R}^N$. By the nonpositivity of $\partial v/\partial n$ we then get

$$\tag{2.4} \int_{\partial\Omega \cap B_r(x)} \left|\frac{\partial v}{\partial n}\right|\,\mathrm{d}\sigma \leq -\int_{\partial\Omega} \psi \frac{\partial v}{\partial n}\,\mathrm{d}\sigma \leq |\mu|(B_{2r}(x)) + C \int_{B_{2r}(x)} \frac{v}{d^2_{\partial\Omega}}.$$

The set

$$E_1 := \left\{x \in \mathbb{R}^N : \limsup_{r \to 0} \frac{1}{r^{N-1}} \int_{B_r(x)} \frac{v}{d^2_{\partial\Omega}} > 0\right\}$$

satisfies $\mathcal{H}^{N-1}(E_1) = 0$, where $\mathcal{H}^{N-1}$ is the Hausdorff measure of $E_1$ of dimension $N - 1$; see e.g. [27, Theorem 2.10]. Since $\mu \equiv 0$ on $\partial\Omega$, by outer regularity of $\mu$ the set

$$E_2 := \left\{x \in \mathbb{R}^N : \limsup_{r \to 0} \frac{|\mu|(B_r(x))}{r^{N-1}} > 0\right\}$$

also satisfies $\mathcal{H}^{N-1}(E_2) = 0$, with the same proof as for $E_1$. Dividing both sides of (2.4) by $r^{N-1}$, it follows that for every $x \in \partial\Omega \setminus (E_1 \cup E_2)$ we have

$$\lim_{r \to 0} \frac{1}{r^{N-1}} \int_{\partial\Omega \cap B_r(x)} \left|\frac{\partial v}{\partial n}\right| d\sigma = 0,$$

and then $\partial v/\partial n = 0$ almost everywhere on $\partial\Omega$ as claimed. □

The choice of the Sobolev space $W_0^{1,1}(\Omega)$ to define the normal derivative is sufficient for our purposes, but one might be interested in a condition that does not require the weak (distributional) derivative to be in $L^1(\Omega; \mathbb{R}^N)$. In fact, the presentation above easily adapts to functions $u \in L^1(\Omega)$ which vanish on the boundary in the sense that

(2.5) $$\lim_{r \to 0} \frac{1}{r} \int_{\{x \in \Omega : d_{\partial\Omega}(x) < r\}} |u| = 0.$$

The reason is that any function $v \in L^1(\Omega)$ such that (2.5) holds and $\Delta v$ is a finite measure in $\Omega$ necessarily belongs to $W_0^{1,1}(\Omega)$; see [42, Propositions 6.3 and 20.2]. Therefore, a family $\widetilde{\mathcal{G}}_u$ defined in terms of (2.5) coincides with our class $\mathcal{G}_u$. An interesting aspect of (2.5) is that such a condition automatically holds for any function $u \in L^1(\Omega)$ that satisfies (2.3).

3. THE DIRICHLET PROBLEM WITH NONHOMOGENEOUS DATA

Given $f \in L^\infty(\Omega)$ and $g \in L^\infty(\partial\Omega)$, the concept of solution of the nonhomogeneous Dirichlet problem for the Schrödinger operator $-\Delta + V$ with $V \in L^1_{\text{loc}}(\Omega)$:

(3.1) $$\begin{cases} -\Delta v + Vv = f & \text{in } \Omega, \\ v = g & \text{on } \partial\Omega, \end{cases}$$

can be straightforwardly defined by $L^1$–$L^\infty$ duality using as test function the solution of the Dirichlet problem

(3.2) $$\begin{cases} -\Delta \zeta + V\zeta = \mu & \text{in } \Omega, \\ \zeta = 0 & \text{on } \partial\Omega, \end{cases}$$

involving $\mu \in L^1(\Omega)$, in the spirit of Stampacchia's definition of weak solutions [45]. In this case, a solution $v \in L^\infty(\Omega)$ of (3.1) is meant to satisfy the identity

(3.3) $$\int_\Omega v\mu = \int_\Omega f\zeta - \int_{\partial\Omega} g \frac{\partial \zeta}{\partial n} d\sigma \quad \text{for every } \mu \in L^1(\Omega).$$

While this notion is enough to investigate the Hopf lemma involving smooth supersolutions of $-\Delta + V$, to deal with nonsmooth ones we rely on a larger class of test functions. Namely, we allow any solution of (3.2) involving finite measures $\mu$ in $\Omega$ which are diffuse with respect to the $W^{1,2}$ capacity. The main result of this section ensures the existence of solutions of (3.1) in this stronger sense:

**Proposition 3.1.** *Let $V \in L^1_{\text{loc}}(\Omega)$ be a nonnegative function. Given $f \in L^\infty(\Omega)$ and $g \in L^\infty(\partial\Omega)$, there exists $v \in W^{1,2}_{\text{loc}}(\Omega) \cap L^\infty(\Omega)$ such that*

$$\text{(3.4)} \qquad \int_\Omega \widehat{v}\,\mathrm{d}\mu = \int_\Omega f\zeta - \int_{\partial\Omega} g\frac{\partial\zeta}{\partial n}\,\mathrm{d}\sigma,$$

*for every finite measure $\mu$ in $\Omega$ which is diffuse with respect to the $W^{1,2}$ capacity, where $\widehat{v}$ is the precise representative of $v$ and $\zeta \in W^{1,1}_0(\Omega)$ satisfies*

$$\text{(3.5)} \qquad -\Delta\zeta + V\zeta = \mu \quad \text{in the sense of distributions in } \Omega.$$

*In particular, there exists a constant $C > 0$ such that*

$$\text{(3.6)} \qquad \|v\|_{L^\infty(\Omega)} \le C\big(\|f\|_{L^\infty(\Omega)} + \|g\|_{L^\infty(\partial\Omega)}\big).$$

We recall that the $W^{1,2}$ capacity of a compact subset $K \subset \mathbb{R}^N$ is defined as

$$\text{cap}_{W^{1,2}}(K) = \inf\big\{\|\varphi\|^2_{W^{1,2}(\mathbb{R}^N)} : \varphi \in C^\infty_c(\mathbb{R}^N) \text{ is nonnegative and } \varphi > 1 \text{ on } K\big\}.$$

A Borel measure $\mu$ in $\Omega$ is diffuse with respect to the $W^{1,2}$ capacity whenever $|\mu|(K) = 0$ for every compact subset $K$ such that $\text{cap}_{W^{1,2}}(K) = 0$. This is the analogue of the notion of absolute continuity between two measures from Measure theory.

We say that $x \in \Omega$ is a Lebesgue point of $v$ and $\widehat{v}(x) \in \mathbb{R}$ is the value of the precise representative of $v$ at $x$ whenever

$$\lim_{r \to 0} \fint_{B_r(x)} |v - \widehat{v}(x)| = 0,$$

where $\fint_{B_r(x)}$ denotes the average integral over $B_r(x)$. For an $L^1_{\text{loc}}$ function, the exceptional set (i.e. the complement of the Lebesgue set in $\Omega$) has Lebesgue measure zero. Since in our case $v$ is a $W^{1,2}_{\text{loc}}$ function, the exceptional set of $v$ is typically smaller and has $W^{1,2}$ capacity zero; see [27, Theorem 4.19; 42, Proposition 8.6]. Thus, the exceptional set of $v$ is irrelevant for diffuse measures.

The existence of solutions of the Dirichlet problem (3.2) for finite diffuse measures $\mu$ is proved in [37] for potentials $V \in L^1_{\text{loc}}(\Omega)$ and depends upon a contraction estimate, following an idea of Brezis and Strauss [16]:

**Proposition 3.2.** *Let $V \in L^1_{\text{loc}}(\Omega)$ be a nonnegative function. For every finite measure $\mu$ in $\Omega$ which is diffuse with respect to the $W^{1,2}$ capacity, there exists a unique function $\zeta \in W^{1,1}_0(\Omega)$ that satisfies equation (3.5). In addition, $V\zeta \in L^1(\Omega)$ and*

$$\text{(3.7)} \qquad \|V\zeta\|_{L^1(\Omega)} \le \|\mu\|_{\mathcal{M}(\Omega)}.$$

The contraction estimate (3.7) holds for any solution of the Dirichlet problem (3.5) provided that $V$ is nonnegative, independently of the fact that $\mu$ is diffuse or not; see e.g. [13, Proposition 4.B.3]. This is a formal consequence of using $\text{sgn}\,\zeta$ as a test function.

The classical weak maximum principle implies by duality the estimate $\|\zeta\|_{L^1(\Omega)} \leq C'\|\Delta\zeta\|_{\mathcal{M}(\Omega)}$. Thus, as a consequence of (3.7) and the identity (3.5) satisfied by $\Delta\zeta$,

$$(3.8) \quad \frac{1}{C'}\|\zeta\|_{L^1(\Omega)} \leq \|\Delta\zeta\|_{\mathcal{M}(\Omega)} \leq \|\mu\|_{\mathcal{M}(\Omega)} + \|V\zeta\|_{L^1(\Omega)} \leq 2\|\mu\|_{\mathcal{M}(\Omega)}.$$

The existence of the distributional normal derivative $\partial\zeta/\partial n$ in $L^1(\partial\Omega)$ relies on the estimate $\|\partial\zeta/\partial n\|_{L^1(\partial\Omega)} \leq \|\Delta\zeta\|_{\mathcal{M}(\Omega)}$. Proceeding as in (3.8) we also get

$$(3.9) \quad \left\|\frac{\partial\zeta}{\partial n}\right\|_{L^1(\partial\Omega)} \leq \|\Delta\zeta\|_{\mathcal{M}(\Omega)} \leq 2\|\mu\|_{\mathcal{M}(\Omega)}.$$

To understand the role played by diffuse measures and the $W^{1,2}$ capacity in this problem, one should keep in mind a classical result in Potential theory which says that for every such a measure $\mu$ one can find a sequence $(\mu_k)_{k\in\mathbb{N}}$ of finite measures that converges strongly to $\mu$ in the space of finite Borel measures $\mathcal{M}(\Omega)$ and such that the solution of the Dirichlet problem

$$(3.10) \quad \begin{cases} -\Delta w_k = \mu_k & \text{in } \Omega, \\ w_k = 0 & \text{on } \partial\Omega, \end{cases}$$

is a bounded function for every $k \in \mathbb{N}$. Those measures can be obtained for example as an application of the Hahn–Banach theorem in the spirit of [17, 28]; see [43, Proposition 2.1] for details. Another strategy relies on the Frostman–Maria boundedness principle by taking $\mu_k := \mu\lfloor_{E_k}$, where $E_k$ is a sublevel set of the solution of the Dirichlet problem (3.10) with datum $\mu$; see Lemma 13.2 in [42] and the remark following that statement. The property that $w_k \in L^\infty(\Omega)$ and $\Delta w_k \in \mathcal{M}(\Omega)$ implies by interpolation that $w_k \in W_0^{1,2}(\Omega)$, hence $\mu_k$ acts as an element in the dual space $(W_0^{1,2}(\Omega))'$.

The existence of a variational solution of the Dirichlet problem (3.2) with better datum $\mu \in (W_0^{1,2}(\Omega))'$ relies on a standard variational approach based on the minimization of the functional

$$(3.11) \quad E(z) = \frac{1}{2}\int_\Omega (|\nabla z|^2 + Vz^2) - \mu[z]$$

in the class $W_0^{1,2}(\Omega) \cap L^2(\Omega; V\,\mathrm{d}x)$. The unique minimizer $\zeta$ satisfies the Euler-Lagrange equation

$$(3.12) \quad \int_\Omega (\nabla\zeta \cdot \nabla z + V\zeta z) = \mu[z],$$

for every $z \in W_0^{1,2}(\Omega) \cap L^2(\Omega; V\,\mathrm{d}x)$. Since $V \in L^1_{\mathrm{loc}}(\Omega)$, the set $C_c^\infty(\Omega)$ is contained in the minimization class $W_0^{1,2}(\Omega) \cap L^2(\Omega; V\,\mathrm{d}x)$. The Euler-Lagrange equation implies in this case that

$$(3.13) \quad -\Delta\zeta + V\zeta = \mu \quad \text{in the sense of distributions in } \Omega.$$

When $\mu$ is in addition a finite measure in $\Omega$, one deduces that $V\zeta \in L^1(\Omega)$ using the test function $z = T_\epsilon(\zeta)/\epsilon$, where $T_\epsilon : \mathbb{R} \to \mathbb{R}$ is the truncation function at $\pm\epsilon$ defined for $t \in \mathbb{R}$ by

$$T_\epsilon(t) = \begin{cases} -\epsilon & \text{if } t < -\epsilon, \\ t & \text{if } -\epsilon \leq t \leq \epsilon, \\ \epsilon & \text{if } t > \epsilon. \end{cases}$$

Indeed, $z$ satisfies $\nabla\zeta \cdot \nabla z \geq 0$ and $|z| \leq 1$. Applying the Euler-Lagrange equation (3.12) with $z$ as above, and letting $\epsilon$ tend to zero, the contraction estimate (3.7) follows from Fatou's lemma. Equation (3.13) thus implies that $\Delta\zeta$ is also a finite measure in $\Omega$. One may prove that $\Delta\zeta$ is a finite measure even for nonnegative Borel functions $V$, although equation (3.13) need not be satisfied in this case; see Section 8 below.

We illustrate these tools with a sketch of the proof of Proposition 3.2:

*Proof of Proposition 3.2.* By linearity, it suffices to consider the case where $\mu$ is nonnegative. Take a sequence of measures $(\mu_k)_{k\in\mathbb{N}}$ converging strongly to $\mu$ in $\mathcal{M}(\Omega)$ such that the solution $w_k$ of the Dirichlet problem (3.10) with density $\mu_k$ is bounded. Hence, $\mu_k$ belongs to $(W_0^{1,2}(\Omega))'$ and then the Dirichlet problem (3.2) with datum $\mu_k$ has a unique solution $\zeta_k \in W_0^{1,2}(\Omega)$. By linearity of the equation (3.13), we have the contraction estimate

$$\|V\zeta_k - V\zeta_l\|_{L^1(\Omega)} \leq \|\mu_k - \mu_l\|_{\mathcal{M}(\Omega)},$$

for every $k, l \in \mathbb{N}$. Hence, the sequence $(V\zeta_k)_{k\in\mathbb{N}}$ is Cauchy in $L^1(\Omega)$. Therefore, $(\Delta\zeta_k)_{k\in\mathbb{N}}$ converges strongly in $\mathcal{M}(\Omega)$, hence by the elliptic estimates of Littman, Stampacchia and Weinberger [34; 42, Proposition 5.1] the sequence $(\zeta_k)_{k\in\mathbb{N}}$ converges strongly in $W_0^{1,p}(\Omega)$ for every $1 \leq p < \frac{N}{N-1}$. In particular, its limit $\zeta$ belongs to $W_0^{1,1}(\Omega)$ and satisfies equation (3.13). $\square$

Before proving Proposition 3.1, we first address the question of existence of solutions that includes the easier $L^1$–$L^\infty$ duality setting and we also develop an approximation scheme of solutions that will be used in the next section:

**Lemma 3.3.** *Let $(g_k)_{k\in\mathbb{N}}$ be a uniformly bounded sequence in $C^\infty(\partial\Omega)$ that converges almost everywhere to $g \in L^\infty(\partial\Omega)$. Then, for every $f \in L^\infty(\Omega)$ and $k \in \mathbb{N}$, there exist $f_k \in L^1_{\text{loc}}(\Omega)$ and $v_k \in W^{1,2}(\Omega) \cap L^\infty(\Omega)$ such that, for every minimizer $\zeta \in W_0^{1,2}(\Omega) \cap L^2(\Omega; V\,\mathrm{d}x)$ of the energy functional (3.11) with datum $\mu \in (W_0^{1,2}(\Omega))' \cap \mathcal{M}(\Omega)$, we have*

*(i) $f_k \in L^1(\Omega; |\zeta|\,\mathrm{d}x)$, $v_k = g_k$ in the sense of traces on $\partial\Omega$ and*

$$\int_\Omega (\nabla v_k \cdot \nabla z + V v_k z) = \int_\Omega f_k z \quad \text{for every } z \in W_0^{1,2}(\Omega) \cap L^1(\Omega; V\,\mathrm{d}x),$$

*(ii) $(f_k)_{k\in\mathbb{N}}$ converges to $f$ in $L^1(\Omega; |\zeta|\,\mathrm{d}x)$,*

$(iii)$ $(v_k)_{k\in\mathbb{N}}$ is uniformly bounded and converges in $L^1(\Omega)$ to the $L^1$–$L^\infty$ duality solution of (3.1).

*Proof of Lemma 3.3.* We construct the function $v_k$ of the form $v_k = u_k + \psi_k$ where $\psi_k \in C^\infty(\overline{\Omega})$ is the harmonic extension of $g_k$ to $\overline{\Omega}$ and $u_k \in W_0^{1,2}(\Omega)$ satisfies

$$-\Delta u_k + V u_k = f - T_k(V)\psi_k \quad \text{in the sense of distributions in } \Omega.$$

Our motivation is that $v_k$ formally satisfies

$$\begin{cases} -\Delta v_k + V v_k = f + (V - T_k(V))\psi_k & \text{in } \Omega, \\ v_k = g_k & \text{on } \partial\Omega, \end{cases}$$

with a warning concerning the fact that $V\psi_k$ need not belong to $L^1(\Omega)$.

Let

$$(3.14) \qquad f_k := f + (V - T_k(V))\psi_k \in L^1_{\text{loc}}(\Omega).$$

The assumption $\mu \in \mathcal{M}(\Omega)$ implies that $V\zeta \in L^1(\Omega)$; see (3.8). Since both $f$ and $\psi_k$ are bounded, we then have $f_k \in L^1(\Omega; |\zeta|\,\mathrm{d}x)$ and the sequence $(f_k)_{k\in\mathbb{N}}$ converges to $f$ in this space by the Dominated convergence theorem. We now show that

$$(3.15) \qquad \int_\Omega (\nabla v_k \cdot \nabla z + V v_k z) = \int_\Omega f_k z \quad \text{for every } z \in W_0^{1,2}(\Omega) \cap L^1(\Omega; V\,\mathrm{d}x),$$

and in particular with $z = \zeta$.

On the one hand, we observe that $u_k$ can be obtained by minimization of the energy functional (3.11) with datum $f - T_k(V)\psi_k$. The Euler-Lagrange equation satisfied by $u_k$ with test function $z$ gives in this case

$$(3.16) \qquad \int_\Omega (\nabla u_k \cdot \nabla z + V u_k z) = \int_\Omega (f - T_k(V)\psi_k) z.$$

On the other hand, since $z \in W_0^{1,2}(\Omega)$ and $\psi_k$ is the harmonic extension of $g_k$,

$$(3.17) \qquad \int_\Omega \nabla \psi_k \cdot \nabla z = -\int_\Omega \Delta \psi_k\, z = 0.$$

For $z \in L^1(\Omega; V\,\mathrm{d}x)$, the integral $\int_\Omega V\psi_k z$ is finite since $\psi_k$ is bounded. Thus adding this integral on both sides of (3.16) and using (3.17) we get (3.15).

We now prove that $(v_k)_{k\in\mathbb{N}}$ is uniformly bounded. To this end, it suffices to establish that $(u_k)_{k\in\mathbb{N}}$ has such a property, and this follows from the pointwise estimate

$$(3.18) \qquad |u_k| \leq \|\psi_k\|_{L^\infty(\Omega)} + \|f\|_{L^\infty(\Omega)}\zeta \quad \text{almost everywhere in } \Omega,$$

where $\zeta \in C_0^\infty(\overline{\Omega})$ is the solution of the Dirichlet problem (1.1) with constant density $\mu \equiv 1$. Indeed, the function $Z_k := u_k - \|\psi_k\|_{L^\infty(\Omega)} - \|f\|_{L^\infty(\Omega)}\zeta$

satisfies

$$\Delta Z_k = V u_k + T_k(V)\psi_k - f + \|f\|_{L^\infty(\Omega)} \quad \text{in the sense of distributions in } \Omega.$$

Thus, by the classical Kato's inequality [31; 42, Proposition 6.6] and non-negativity of $V$,

$$\Delta Z_k^+ \geq \chi_{\{Z_k \geq 0\}}\big(V u_k + T_k(V)\psi_k - f + \|f\|_{L^\infty(\Omega)}\big) \geq 0$$

in the sense of distributions in $\Omega$. Since $Z_k^+ \in W_0^{1,2}(\Omega)$, the weak maximum principle implies that $Z_k^+ \leq 0$ almost everywhere in $\Omega$; see e.g. [42, Propositions 6.1 and 6.5]. Therefore, $Z_k^+$ vanishes in $\Omega$ and we get

$$u_k \leq \|\psi_k\|_{L^\infty(\Omega)} + \|f\|_{L^\infty(\Omega)}\zeta \quad \text{almost everywhere in } \Omega.$$

A similar estimate holds for $-u_k$ and one deduces (3.18). As $\|\psi_k\|_{L^\infty(\Omega)} = \|g_k\|_{L^\infty(\partial\Omega)}$ and the sequence $(g_k)_{k\in\mathbb{N}}$ is uniformly bounded, we deduce that $(v_k)_{k\in\mathbb{N}}$ is uniformly bounded. This type of property where the potential of the Schrödinger operator forces the equation to have bounded solutions from data that are merely $L^1$ has been further investigated by Arcoya and Boccardo [7], based on suitable choices of test functions.

We are left with the convergence of the sequence $(v_k)_{k\in\mathbb{N}}$. We have just proved that $(v_k)_{k\in\mathbb{N}}$ is uniformly bounded. Since $V \in L^1_{\mathrm{loc}}(\Omega)$ and

$$-\Delta v_k + V v_k = f_k \quad \text{in the sense of distributions in } \Omega,$$

the sequence $(\Delta v_k)_{k\in\mathbb{N}}$ is bounded in $L^1(\omega)$ for every $\omega \Subset \Omega$. It then follows by interpolation that $(v_k)_{k\in\mathbb{N}}$ is bounded in $W^{1,2}(\omega)$ for every $\omega \Subset \Omega$. Thus, there exists a subsequence $(v_{k_j})_{j\in\mathbb{N}}$ that converges to some function $v \in W^{1,2}_{\mathrm{loc}}(\Omega)$ weakly in $W^{1,2}(\omega)$ for every $\omega \Subset \Omega$ and strongly in $L^1(\Omega)$; the latter holds in the entire domain $\Omega$ by uniform boundedness of $(v_k)_{k\in\mathbb{N}}$.

To identify the limit $v$, we return to (3.16) and (3.17) to prove that under the additional assumption that $\mu \in L^\infty(\Omega)$, one has

(3.19) $$\int_\Omega v_k \mu = \int_\Omega f_k \zeta - \int_{\partial\Omega} g_k \frac{\partial \zeta}{\partial n}\,\mathrm{d}\sigma.$$

For such a $\mu$, the quantity $\mu[u_k]$ can be computed through integration of $u_k\mu$. From the Euler-Lagrange equation satisfied by $\zeta$ with test function $u_k$ and (3.16), we get

$$\int_\Omega u_k \mu = \mu[u_k] = \int_\Omega (f - T_k(V)\psi_k)\zeta.$$

Since $v_k = g$ on $\partial\Omega$ and $\zeta$ has a distributional normal derivative, (3.17) with $z = \zeta$ implies that

$$-\int_\Omega \psi_k \Delta\zeta + \int_{\partial\Omega} g_k \frac{\partial\zeta}{\partial n}\,\mathrm{d}\sigma = 0.$$

Thus adding $\int_\Omega V\psi_k\zeta$ on both sides and using the fact that $\mu = -\Delta\zeta + V\zeta$ in the sense of measures in $\Omega$, we get

$$\int_\Omega \psi_k \mu = \int_\Omega V\psi_k\zeta - \int_{\partial\Omega} g_k \frac{\partial \zeta}{\partial n}\,\mathrm{d}\sigma.$$

A combination of the first and third identities implies (3.19). As $k = k_j$ tends to infinity, we have that $v$ satisfies

$$\int_\Omega v\mu = \int_\Omega f\zeta - \int_{\partial\Omega} g \frac{\partial\zeta}{\partial n}\,\mathrm{d}\sigma \quad \text{for every } \mu \in L^\infty(\Omega).$$

This already gives the uniqueness of the limit and in particular the entire sequence $(v_k)_{k\in\mathbb{N}}$ converges to $v$ in $L^1(\Omega)$. That this identity holds for every $\mu \in L^1(\Omega)$ follows from approximation of $\mu$ by bounded functions $(\mu_k)_{k\in\mathbb{N}}$. Indeed, the solutions $(\zeta_k)_{k\in\mathbb{N}}$ associated to that sequence converge to $\zeta$ in $L^1(\Omega)$ and $(\partial\zeta_k/\partial n)_{k\in\mathbb{N}}$ converges to $\partial\zeta/\partial n$ in $L^1(\partial\Omega)$ by estimates (3.8) and (3.9). It thus suffices to use $\mu_k$ as test function and let $k$ tend to infinity. The proof of the lemma is complete. $\square$

A finite measure $\nu$ in $\Omega$ that belongs to the dual space $(W_0^{1,2}(\Omega))'$ verifies

(3.20) $$\left| \int_\Omega \varphi\,\mathrm{d}\nu \right| \leq C\|\varphi\|_{W^{1,2}(\Omega)},$$

for every $\varphi \in C_c^\infty(\Omega)$. By density of $C_c^\infty(\Omega)$ in $W_0^{1,2}(\Omega)$, the linear functional

$$\varphi \in C_c^\infty(\Omega) \longmapsto \int_\Omega \varphi\,\mathrm{d}\nu$$

has a unique continuous extension to $W_0^{1,2}(\Omega)$. Denoting such an extension by $\nu[u]$ for every $u \in W_0^{1,2}(\Omega)$, one can represent $\nu[u]$ as integration of $u$ with respect to $\nu$. Indeed, estimate (3.20) implies that $\nu$, as a measure, is diffuse with respect to the $W^{1,2}$ capacity, the precise representative $\widehat{u}$ has an exceptional set with $W^{1,2}$ capacity zero, and $\widehat{u} \in L^1(\Omega;\nu)$; see e.g. [30; 42, Proposition 16.5]. Moreover,

(3.21) $$\nu[u] = \int_\Omega \widehat{u}\,\mathrm{d}\nu \quad \text{for every } u \in W_0^{1,2}(\Omega).$$

*Proof of Proposition 3.1.* Estimate (3.1) is a straightforward consequence of the $L^1$–$L^\infty$ duality. Indeed, for any $\mu \in L^1(\Omega)$, by estimates (3.8) and (3.9) we have

$$\left|\int_\Omega v\mu\right| = |\mu[v]| \leq \|f\|_{L^\infty(\Omega)}\|\zeta\|_{L^1(\Omega)} + \|g\|_{L^\infty(\partial\Omega)}\left\|\frac{\partial\zeta}{\partial n}\right\|_{L^1(\partial\Omega)}$$
$$\leq 2\big(C'\|f\|_{L^\infty(\Omega)} + \|g\|_{L^\infty(\partial\Omega)}\big)\|\mu\|_{L^1(\Omega)}.$$

By $L^1$–$L^\infty$ duality, we deduce that $v \in L^\infty(\Omega)$ and

$$\|v\|_{L^\infty(\Omega)} \leq 2\big(C'\|f\|_{L^\infty(\Omega)} + \|g\|_{L^\infty(\partial\Omega)}\big).$$

We thus have the conclusion with $C := 2\max\{C', 1\}$.

The proof of Lemma 3.3 may be seen as a first step in establishing Proposition 3.1. We follow the notation there: We recall that $(g_k)_{k\in\mathbb{N}}$ is a uniformly bounded sequence in $C^\infty(\partial\Omega)$ that converges almost everywhere to $g$ and $(f_k)_{k\in\mathbb{N}}$ is defined by (3.14) and converges to $f$ in $L^1(\Omega; |\zeta|\,dx)$, where $\zeta$ is the solution of (3.5) with $\mu \in (W_0^{1,2}(\Omega))' \cap \mathcal{M}(\Omega)$. The sequence $(v_k)_{k\in\mathbb{N}}$ defined by $v_k = u_k + \psi_k$ is uniformly bounded and also bounded in $W^{1,2}(\omega)$ for every open subset $\omega \Subset \Omega$.

We now prove that if $\mu \in (W_0^{1,2}(\Omega))' \cap \mathcal{M}(\Omega)$ has *compact support in $\Omega$*, then

$$(3.22) \quad \mu[v_k\varphi] = \int_\Omega f_k \zeta - \int_{\partial\Omega} g_k \frac{\partial\zeta}{\partial n}\,d\sigma,$$

where $\varphi \in C_c^\infty(\Omega)$ is any function such that $\varphi = 1$ on $\operatorname{supp}\mu$. Proceeding as in the case where $\mu$ was assumed to belong to $L^\infty(\Omega)$, we have

$$(3.23) \quad \mu[u_k] = \int_\Omega (f - T_k(V)\psi_k)\zeta$$

and

$$(3.24) \quad \int_\Omega \psi_k\,d\mu = \int_\Omega V\psi_k\zeta - \int_{\partial\Omega} g_k \frac{\partial\zeta}{\partial n}\,d\sigma.$$

For $\varphi \in C_c^\infty(\Omega)$ such that $\varphi = 1$ on $\operatorname{supp}\mu$, we also have

$$(3.25) \quad \mu[u_k] = \mu[u_k\varphi] \quad\text{and}\quad \int_\Omega \psi_k\,d\mu = \int_\Omega \psi_k\varphi\,d\mu = \mu[\psi_k\varphi].$$

A combination of Equations (3.23) to (3.25) then implies (3.22).

Take a subsequence $(v_{k_j})_{j\in\mathbb{N}}$ that converges to $v$ weakly in $W^{1,2}(\omega)$ for every $\omega \Subset \Omega$ and in $L^1(\Omega)$. By (3.21), we have

$$\lim_{j\to\infty} \mu[v_{k_j}\varphi] = \mu[v\varphi] = \int_\Omega \widehat{v\varphi}\,d\mu = \int_\Omega \widehat{v}\,d\mu.$$

In view of the convergences of $(f_k)_{k\in\mathbb{N}}$ and $(g_k)_{k\in\mathbb{N}}$, as $k = k_j$ tends to infinity in (3.22) we conclude that

$$(3.26) \quad \int_\Omega \widehat{v}\,d\mu = \int_\Omega f\zeta - \int_{\partial\Omega} g \frac{\partial\zeta}{\partial n}\,d\sigma,$$

when $\mu \in \mathcal{M}(\Omega) \cap (W_0^{1,2}(\Omega))'$ has compact support in $\Omega$.

We finally prove that $v$ satisfies identity (3.4) for every test function $\zeta$ as in the statement. To this end, we may assume that $\mu$ is nonnegative. As in the proof of Proposition 3.2, take a sequence $(\mu_k)_{k\in\mathbb{N}}$ of finite measures that converges strongly to $\mu$ in $\mathcal{M}(\Omega)$ and such that, for each measure $\mu_k$, the solution $w_k$ of the Dirichlet problem (3.10) with density $\mu_k$ is bounded. We may also assume that each $\mu_k$ has compact support in $\Omega$. By interpolation, $w_k \in W_0^{1,2}(\Omega)$ and then $\mu_k \in (W_0^{1,2}(\Omega))'$. Denoting by $\zeta_k \in W_0^{1,2}(\Omega) \cap L^2(\Omega; V\,dx)$ the solution of (3.5) with $\mu_k$, it follows from (3.26) that

$$\int_\Omega \widehat{v}\,d\mu_k = \int_\Omega f\zeta_k - \int_{\partial\Omega} g \frac{\partial\zeta_k}{\partial n}\,d\sigma.$$

On the one hand, since the function $\widehat{v}$ is bounded, by strong convergence of the sequence $(\mu_k)_{k \in \mathbb{N}}$ we have

$$\lim_{k \to \infty} \int_\Omega \widehat{v} \, \mathrm{d}\mu_k = \int_\Omega \widehat{v} \, \mathrm{d}\mu.$$

On the other hand, by estimates (3.8) and (3.9) the sequence $(\zeta_k)_{k \in \mathbb{N}}$ converges to $\zeta$ in $L^1(\Omega)$ and $(\partial \zeta_k / \partial n)_{k \in \mathbb{N}}$ converges to $\partial \zeta / \partial n$ in $L^1(\partial \Omega)$. By the boundedness of $f$ and $g$ we get

$$\int_\Omega \widehat{v} \, \mathrm{d}\mu = \lim_{k \to \infty} \left( \int_\Omega f \zeta_k - \int_{\partial \Omega} g \frac{\partial \zeta_k}{\partial n} \, \mathrm{d}\sigma \right) = \int_\Omega f \zeta - \int_{\partial \Omega} g \frac{\partial \zeta}{\partial n} \, \mathrm{d}\sigma.$$

The proof is complete. $\square$

## 4. CONSTRUCTION OF POSITIVE TEST FUNCTIONS

Given any nontrivial nonnegative boundary datum $g \in L^\infty(\partial \Omega)$, the main result of this section gives a recipe to construct $f \in L^\infty(\Omega)$ such that $f < 0$ almost everywhere in $\Omega$ while the Dirichlet problem (3.1) with mixed sign datum $(f, g)$ has a *nonnegative* solution.

**Proposition 4.1.** *There exists a bounded continuous function $H : \mathbb{R}_+ \to \mathbb{R}_+$, with $H(t) > 0$ for $t > 0$, such that, for every nonnegative functions $V \in L^1_{\mathrm{loc}}(\Omega)$ and $g \in L^\infty(\partial \Omega)$, if $w$ satisfies the Dirichlet problem*

$$\begin{cases} -\Delta w + Vw = 0 & \text{in } \Omega, \\ w = g & \text{on } \partial \Omega, \end{cases}$$

*and if $v$ satisfies*

$$\begin{cases} -\Delta v + Vv = H(w) & \text{in } \Omega, \\ v = 0 & \text{on } \partial \Omega, \end{cases}$$

*then we have $w \geq v$ almost everywhere in $\Omega$.*

We illustrate this proposition with a proof of Theorem 1 for smooth supersolutions involving potentials in $L^1(\Omega; d_{\partial \Omega} \, \mathrm{d}x)$. An important ingredient is the following strong maximum principle for $L^1$ potentials that was proved independently by Ancona [1] and Trudinger [46]; see also [14]:

**Proposition 4.2.** *Let $V \in L^1_{\mathrm{loc}}(\Omega)$. If $u \in L^1(\Omega)$ is a nonnegative supersolution of the Schrödinger operator $-\Delta + V$ and if $\int_\Omega u > 0$, then $u > 0$ almost everywhere in $\Omega$.*

*Proof of Theorem 1 when $V \in L^1(\Omega; d_{\partial \Omega} \, \mathrm{d}x)$ and $u \in C_0^\infty(\overline{\Omega})$.* Since $u$ is nonnegative, we may assume from the beginning that $V$ is also nonnegative. Assume by contradiction that the set $A := \{\partial u / \partial n = 0\}$ is not negligible with respect to the surface measure on $\partial \Omega$. We solve the Dirichlet problems of Proposition 4.1 starting with the boundary condition $g = \chi_A$.

Using the notation in that statement, the function $w - v$ is nonnegative. Since $w - v$ satisfies the Dirichlet problem with datum $(-H(w), \chi_A)$, using $u$ as test function (i.e. take $u = \zeta$ in (3.3)) we have

$$\int_\Omega (w-v)(-\Delta u + Vu) = \int_\Omega (-H(w))u - \int_{\partial\Omega} \chi_A \frac{\partial u}{\partial n}\,d\sigma.$$

Observe that $u$ is an admissible test function because $V \in L^1(\Omega; d_{\partial\Omega}\,dx)$, which implies that $-\Delta u + Vu \in L^1(\Omega)$. By the choice of $A$ the last integral vanishes. Thus,

$$\int_\Omega H(w)u = -\int_\Omega (w-v)(-\Delta u + Vu) \le 0,$$

by nonnegativity of the integrand in the right-hand side. This implies that $H(w)u = 0$ almost everywhere in $\Omega$. Since $A$ has positive surface measure, $w$ and $v$ are not identically zero; this follows from the fact that one can use as test function in both problems the solution of the Dirichlet problem (1.1) with datum $\mu \equiv 1$ (cf. Proposition 5.1 below). By the strong maximum principle above for the Schrödinger operator $-\Delta + V$, we have that $w > 0$ almost everywhere in $\Omega$, hence $H(w)$ satisfies the same property. Therefore, $u \equiv 0$ and the conclusion follows for smooth supersolutions. $\square$

To prove Proposition 4.1, we need a version of Kato's inequality adapted to the nonhomogeneous Dirichlet problem involving potentials $V \in L^1_{\mathrm{loc}}(\Omega)$.

**Lemma 4.3.** *Let $V \in L^1_{\mathrm{loc}}(\Omega)$ be a nonnegative function. Given $f \in L^\infty(\Omega)$ and $g \in L^\infty(\partial\Omega)$, if $v \in L^\infty(\Omega)$ satisfies the Dirichlet problem (3.1), then*

$$\int_\Omega v^+ \le \int_{\{v>0\}} f\xi - \int_{\partial\Omega} g^+ \frac{\partial \xi}{\partial n}\,d\sigma,$$

*where $\xi \in W_0^{1,2}(\Omega) \cap L^\infty(\Omega)$ is the (nonnegative) solution of the Dirichlet problem (3.2) with datum $\mu \equiv 1$.*

*Proof of Lemma 4.3.* The proof is based on an approximation of $v$ by functions $v_k \in W^{1,2}(\Omega) \cap L^\infty(\Omega)$, following the notation in Lemma 3.3. By the contraction estimate (3.7), $\Delta\xi \in L^1(\Omega)$ and $\xi$ has a distributional normal derivative $\partial\xi/\partial n \in L^1(\partial\Omega)$. Thus,

$$\int_\Omega (\nabla\xi \cdot \nabla\psi + V\xi\psi) = \int_\Omega \psi + \int_{\partial\Omega} \psi \frac{\partial\xi}{\partial n}\,d\sigma \quad \text{for every } \psi \in C^\infty(\overline{\Omega}).$$

Since $\xi \in W_0^{1,2}(\Omega)$ and $V\xi \in L^1(\Omega)$, this identity holds by approximation of $\psi$ for every $\psi \in W^{1,2}(\Omega) \cap L^\infty(\Omega)$. In particular, we can take $\psi = J(v_k)$, where $J : \mathbb{R} \to \mathbb{R}$ is a smooth function to be chosen later on. Since $J(v_k) = J(g_k)$ in the sense of traces on $\partial\Omega$, where $g_k \in C^\infty(\partial\Omega)$ is an approximation of $g$, we get

$$(4.1) \quad \int_\Omega \bigl(J'(v_k)\nabla\xi \cdot \nabla v_k + V\xi J(v_k)\bigr) = \int_\Omega J(v_k) + \int_{\partial\Omega} J(g_k) \frac{\partial\xi}{\partial n}\,d\sigma.$$

On the other hand, applying the Euler-Lagrange equation in the statement of Lemma 3.3 with test function $z = J'(v_k)\xi$ (which satisfies $|z| \le C|\xi|$ and thus $z \in L^1(\Omega; V\,\mathrm{d}x)$), we have

$$\int_\Omega \bigl(J''(v_k)|\nabla v_k|^2 \xi + J'(v_k)\nabla\xi \cdot \nabla v_k + V v_k J'(v_k)\xi\bigr) = \int_\Omega f_k J'(v_k)\xi.$$

Assuming that $J'' \ge 0$, by nonnegativity of $\xi$ we get

(4.2) $$\int_\Omega \bigl(J'(v_k)\nabla v_k \cdot \nabla\xi + V v_k J'(v_k)\xi\bigr) \le \int_\Omega f_k J'(v_k)\xi.$$

Subtracting (4.2) from (4.1), we get

$$\int_\Omega V\xi[J(v_k) - v_k J'(v_k)] \ge \int_\Omega J(v_k) - \int_\Omega f_k J'(v_k)\xi + \int_{\partial\Omega} J(g_k)\frac{\partial\xi}{\partial n}\,\mathrm{d}\sigma.$$

We now take $J$ convex such that $J(t) = 0$ for $t \le 0$ and $0 \le J(t) \le t$ for $t \ge 0$. In particular, for every $t \in \mathbb{R}$ we have $J(t) \le J'(t)t$. Since $V$ and $\xi$ are nonnegative, the integrand in the left-hand side is nonpositive and we deduce that

$$\int_\Omega J(v_k) \le \int_\Omega f_k J'(v_k)\xi - \int_{\partial\Omega} J(g_k)\frac{\partial\xi}{\partial n}\,\mathrm{d}\sigma.$$

The sequence $(f_k)_{k\in\mathbb{N}}$ provided by Lemma 3.3 converges to $f$ in $L^1(\Omega; \xi\,\mathrm{d}x)$. As $k$ tends to infinity, by pointwise convergence and boundedness of $(v_k)_{k\in\mathbb{N}}$ and $(g_k)_{k\in\mathbb{N}}$ and by the Dominated convergence theorem we thus get

$$\int_\Omega J(v) \le \int_\Omega f J'(v)\xi - \int_{\partial\Omega} J(g)\frac{\partial\xi}{\partial n}\,\mathrm{d}\sigma.$$

To conclude, we apply this inequality to a sequence $(J_i)_{i\in\mathbb{N}}$ of convex functions as above that converges pointwise to $t \mapsto t^+$ and such that $(J'_i)_{i\in\mathbb{N}}$ converges pointwise to $\chi_{(0,\infty)}$. As $i$ tends to infinity, we have the conclusion. $\square$

*Proof of Proposition 4.1.* For every $\epsilon > 0$, we claim that if $z_\epsilon$ satisfies the Dirichlet problem

$$\begin{cases} -\Delta z_\epsilon + V z_\epsilon = \chi_{\{w > \epsilon C\}} & \text{in } \Omega, \\ z_\epsilon = 0 & \text{on } \partial\Omega, \end{cases}$$

for $C > 0$ sufficiently large, then we have $w \ge \epsilon z_\epsilon$. Indeed, by Kato's inequality above applied to $v = \epsilon z_\epsilon - w$, we have

$$\int_\Omega (\epsilon z_\epsilon - w)^+ \le \int_{\{\epsilon z_\epsilon > w\}} \epsilon \chi_{\{w>\epsilon C\}}\xi - \int_{\partial\Omega}(-g)^+ \frac{\partial\xi}{\partial n}\,\mathrm{d}\sigma = \epsilon \int_{\{\epsilon z_\epsilon > w > \epsilon C\}} \xi,$$

since $g$ is assumed to be nonnegative. Observe that

$$z_\epsilon \le \xi \le \|\xi\|_{L^\infty(\Omega)} \quad \text{for every } \epsilon > 0.$$

Taking $C := \|\xi\|_{L^\infty(\Omega)}$, the set $\{\epsilon z_\epsilon > w > \epsilon C\}$ is then negligible and we deduce that
$$\int_\Omega (\epsilon z_\epsilon - w)^+ \leq 0.$$
Hence, $w \geq \epsilon z_\epsilon$ almost everywhere in $\Omega$. Applying this conclusion with $\epsilon = 1/2^k$ for every $k \in \mathbb{N}_*$, we get
$$w = \sum_{k=1}^\infty \frac{1}{2^k} w \geq \sum_{k=1}^\infty \frac{1}{2^k} \cdot \frac{1}{2^k} z_{1/2^k} =: \widetilde{v},$$
where $\widetilde{v}$ satisfies the Dirichlet problem
$$\begin{cases} -\Delta \widetilde{v} + V \widetilde{v} = \widetilde{f} & \text{in } \Omega, \\ \widetilde{v} = 0 & \text{on } \partial\Omega, \end{cases}$$
with
$$\widetilde{f}(x) := \sum_{k=1}^\infty \frac{1}{2^k} \cdot \frac{1}{2^k} \chi_{\{w > C/2^k\}}(x) = \sum_{k=1}^\infty \frac{1}{2^k} \cdot \frac{1}{2^k} \chi_{(C/2^k, \infty)}(w(x)).$$
Observe that for every $s \geq 0$,
$$\begin{aligned}\sum_{k=1}^\infty \frac{1}{2^k} \cdot \frac{1}{2^k} \chi_{(C/2^k, \infty)}(s) &\geq \sum_{k=1}^\infty \int_{\frac{1}{2^k}}^{\frac{1}{2^{k-1}}} \frac{t}{2} \chi_{(Ct, \infty)}(s) \, \mathrm{d}t \\ &= \int_0^1 \frac{t}{2} \chi_{(0, s/C)}(t) \, \mathrm{d}t \\ &= \int_0^{\min\{1, s/C\}} \frac{t}{2} \, \mathrm{d}t = \frac{1}{4}\left(\min\{1, s/C\}\right)^2 =: H(s).\end{aligned}$$
By this computation, we thus have $\widetilde{f} \geq H(w)$ in $\Omega$. Hence by comparison it follows that the solution $v$ of the Dirichlet problem with datum $H(w)$ satisfies $\widetilde{v} \geq v$ in $\Omega$. Therefore, $w \geq v$ as we wanted to prove. $\square$

## 5. Nontrivial solutions for the nonhomogeneous problem

The existence of nontrivial solutions of the boundary value problem

(5.1) $$\begin{cases} -\Delta w + Vw = 0 & \text{in } \Omega, \\ w = g & \text{on } \partial\Omega, \end{cases}$$

in the sense of Proposition 3.1 is related to the existence of supersolutions of the Schrödinger operator $-\Delta + V$ with negative normal derivative through the following

**Proposition 5.1.** *Let $V \in L^1_{\mathrm{loc}}(\Omega)$ and $g \in L^\infty(\partial\Omega)$ be nonnegative functions. Then, the (nonnegative) solution $w \in W^{1,2}_{\mathrm{loc}}(\Omega) \cap L^\infty(\Omega)$ of the Dirichlet problem (5.1) with datum $g$ satisfies*
$$\int_\Omega w > 0$$

*if and only if there exists a (nonnegative) supersolution $\zeta \in W_0^{1,1}(\Omega)$ of the Schrödinger operator $-\Delta + V$ such that the measure $-\Delta \zeta + V\zeta$ is finite and diffuse in $\Omega$ and*

$$\int_{\partial \Omega} g \frac{\partial \zeta}{\partial n} \, \mathrm{d}\sigma < 0. \tag{5.2}$$

Proposition 5.1 is used in the proofs of Theorems 1 and 2. We deduce a posteriori from Theorem 2 that once condition (5.2) holds for *one* supersolution, then it holds for *every* nontrivial supersolution $\zeta \in W_0^{1,1}(\Omega)$ such that $V\zeta \in L^1(\Omega)$.

The nonnegativity of $w$ and $\zeta$ follows from Proposition 3.1 and Lemma 4.3. Indeed, by Lemma 4.3 applied to $v := -w$ we have

$$\int_\Omega (-w)^+ \leq -\int_{\partial \Omega} (-g)^+ \frac{\partial \xi}{\partial n} \, \mathrm{d}\sigma = 0.$$

Thus, $-w \leq 0$ almost everywhere in $\Omega$. The same argument applies to solutions of (3.1) such that $f \geq 0$ in $\Omega$ and $g \geq 0$ on $\partial \Omega$. In particular, by Proposition 3.1 and the nonnegativity of the solutions of (3.1) in this case we have

$$\int_\Omega f\zeta - \int_{\partial \Omega} g \frac{\partial \zeta}{\partial n} \, \mathrm{d}\sigma = \int_\Omega \widehat{v} \, \mathrm{d}\mu \geq 0,$$

where $\mu = -\Delta \zeta + V\zeta$. Since this is true for every nonnegative $f$ and $g$, we deduce that $\zeta \geq 0$ in $\Omega$ and $\partial \zeta / \partial n \leq 0$ on $\partial \Omega$.

*Proof of Proposition 5.1.* "$\Longleftarrow$" Since $\partial \zeta / \partial n \leq 0$ almost everywhere on $\partial \Omega$, by Proposition 3.1 we have

$$\int_\Omega \widehat{w} \, \mathrm{d}\mu = -\int_{\partial \Omega} g \frac{\partial \zeta}{\partial n} \, \mathrm{d}\sigma > 0.$$

In particular, $w$ is a nonzero solution of (5.1). Since $g$ is nonnegative on $\partial \Omega$, $w$ is nonnegative in $\Omega$ and the implication follows.

"$\Longrightarrow$" Since $w$ is a nontrivial nonnegative solution, there exists a compact subset $K \subset \Omega$ with positive $W^{1,2}$ capacity which is contained in the Lebesgue set of $w$ and is such that $\widehat{w} \geq \epsilon$ on $K$ for some constant $\epsilon > 0$. Take a finite nonnegative diffuse measure $\mu$ supported on $K$ such that $\mu(K) > 0$. The existence of such a measure follows from the Hahn-Banach theorem; see e.g. [42, Proposition A.17]. We take as supersolution the function $\zeta \in W_0^{1,1}(\Omega)$ such that

$$-\Delta \zeta + V\zeta = \mu \quad \text{in the sense of distributions in } \Omega;$$

see Proposition 3.2. Applying $\zeta$ as a test function in the Dirichlet problem satisfied by $w$, we have

$$\epsilon \mu(K) \leq \int_\Omega \widehat{w} \, \mathrm{d}\mu = -\int_{\partial \Omega} g \frac{\partial \zeta}{\partial n} \, \mathrm{d}\sigma,$$

which concludes the proof. $\square$

The previous proposition raises the question of how to construct supersolutions of the Schrödinger operator $-\Delta + V$ with pointwise control on its distributional normal derivative to ensure that (5.2) is satisfied.

**Proposition 5.2.** *Let $V \in L^1_{\mathrm{loc}}(\Omega)$ be a nonnegative function. For every $w \in W^{1,1}_0(\Omega)$ such that $\Delta w$ is a finite measure and $Vw \in L^1(\Omega)$, there exists a nonnegative function $\widetilde{w} \in W^{1,2}_0(\Omega) \cap L^\infty(\Omega)$ such that*

*($O_1$) $\Delta \widetilde{w}$ is a finite measure in $\Omega$ and $V\widetilde{w} \in L^1(\Omega)$,*
*($O_2$) $\partial \widetilde{w}/\partial n \leq \partial w/\partial n$ almost everywhere on $\partial \Omega$,*
*($O_3$) $-\Delta \widetilde{w} + V\widetilde{w} \geq 0$ in the sense of distributions in $\Omega$.*

For example, when $V$ is a nonnegative Hopf potential, then taking $w := v$ as the function given by Proposition 2.5 with $u := \zeta_0$ one finds, as an application of Proposition 5.2, a supersolution of $-\Delta + V$ having a distributional normal derivative $\partial \widetilde{w}/\partial n < 0$ almost everywhere on $\partial \Omega$.

*Proof of Proposition 5.2.* Since $w \in W^{1,1}_0(\Omega)$ and the measure $\Delta w$ is finite in $\Omega$, by interpolation we have that $T_1(w) \in W^{1,2}_0(\Omega)$. Since $Vw \in L^1(\Omega)$, we also have $VT_1(w) \in L^1(\Omega)$. By Kato's inequality up to the boundary (2.1), $\Delta T_1(w)$ is a finite measure in $\Omega$. By (2.2), we also have

$$\text{(5.3)} \qquad \frac{\partial T_1(w)}{\partial n} = \frac{\partial w}{\partial n}.$$

Thus, replacing $w$ by $T_1(w)$ if necessary, we may henceforth assume that $w \in W^{1,2}_0(\Omega) \cap L^\infty(\Omega)$. The measure $\Delta w$ in this case is diffuse with respect to the $W^{1,2}$ capacity, hence the measure

$$\nu := -\Delta w + Vw$$

is also finite and diffuse in $\Omega$, and so is its positive part $\nu^+$.

By Proposition 3.2, the Dirichlet problem with nonnegative potential $V$:

$$\begin{cases} -\Delta z + Vz = \nu^+ & \text{in } \Omega, \\ z = 0 & \text{on } \partial \Omega, \end{cases}$$

has a solution. It is not clear for example why $z$ is bounded, for this reason we now prove that $\widetilde{w} := T_1(z)$ satisfies the required properties. The contraction estimate (3.7) implies that $Vz \in L^1(\Omega)$, hence the measure $\Delta z$ is finite. Proceeding as in the first part of the proof, we have $\widetilde{w} \in W^{1,2}_0(\Omega)$, $\Delta \widetilde{w}$ is a finite measure in $\Omega$, and

$$\frac{\partial \widetilde{w}}{\partial n} = \frac{\partial z}{\partial n}.$$

By the comparison principle between solutions of the Dirichlet problem we have $z \geq w$ in $\Omega$. Then, by comparison between normal derivatives,

$$\frac{\partial \widetilde{w}}{\partial n} = \frac{\partial z}{\partial n} \leq \frac{\partial w}{\partial n},$$

which is $(O_2)$. Since $z$ is nonnegative and $\Delta z \leq Vz$ in the sense of distributions in $\Omega$, a straightforward variant of Kato's inequality yields

$$\Delta \widetilde{w} = \Delta(\min\{z,1\}) \leq \chi_{\{z<1\}} Vz \tag{5.4}$$

in the sense of distributions in $\Omega$; see [42, Proposition 6.9]. By nonnegativity of $V$, the right-hand side is smaller than $V\widetilde{w}$, and we deduce that $\widetilde{w}$ satisfies $(O_3)$. The proof is complete. □

Ancona's argument leading to (2.2) is based on tools from Potential theory. There is another strategy which allows one to prove a smooth counterpart of this formula based on a PDE approach, which is enough to prove Proposition 5.2. More precisely, given a smooth function $\Phi : \mathbb{R} \to \mathbb{R}$ such that $\Phi''$ has compact support, it has been proved in [20] using the notion of renormalized solution that, for every $u \in W_0^{1,1}(\Omega)$ such that $\Delta u$ is a finite measure in $\Omega$, one has that $\Delta \Phi(u)$ is also a finite measure in $\Omega$ and the following holds:

$$\Delta \Phi(u) = \Phi'(u)(\Delta u)_{\mathrm{d}} + \Phi''(u)|\nabla u|^2. \tag{5.5}$$

Here, $(\Delta u)_{\mathrm{d}}$ denotes the part of the measure $\Delta u$ which is diffuse with respect to the $W^{1,2}$ capacity that arises from the Lebesgue decomposition of measures.

The approximation scheme from [15] to prove that the distributional normal derivative belongs to $L^1(\partial\Omega)$ is based on a strong approximation of the measure $\Delta u$ by measures with compact support. In this case, one deduces using the identity (5.5) that the solutions $u_k$ of the approximating Dirichlet problems are such that $(\Delta \Phi(u_k))_{k \in \mathbb{N}}$ converges strongly to $\Delta \Phi(u)$ in $\Omega$, and one then deduces that

$$\frac{\partial \Phi(u)}{\partial n} = \Phi'(0) \frac{\partial u}{\partial n}.$$

In particular, if $\Phi$ is an approximation of the truncation function such that $\Phi'(0) = 1$, one gets an equality between the normal derivatives as in (2.2) and (5.3).

## 6. Proofs of the main results

*Proof of Theorem 2.* Since $u_1$ and $u_2$ are nonnegative, they are also supersolutions for the Schrödinger operator $-\Delta + V^+$ with nonnegative potential. We may thus assume from the beginning that $V$ is nonnegative. We split the proof into three steps:

*Step* 1. We prove the theorem under the additional assumption that the measures $\Delta u_1$ and $\Delta u_2$ are finite and diffuse.

By Proposition 2.1, both $\partial u_1/\partial n$ and $\partial u_2/\partial n$ are nonpositive on $\partial\Omega$. Let us prove that

(6.1) $$\frac{\partial u_1}{\partial n} < 0 \quad \text{almost everywhere on } \left\{\frac{\partial u_2}{\partial n} < 0\right\}$$

with respect to the surface measure on $\partial\Omega$. Assume by contradiction that there exists a Borel set $A \subset \{\partial u_1/\partial n = 0\}$ such that

$$\int_A \frac{\partial u_2}{\partial n}\,d\sigma < 0.$$

By Proposition 5.1, the solution $w$ of the Dirichlet problem (5.1) with datum $g = \chi_A$ is nontrivial. Since $\Omega$ is connected, by the strong maximum principle for $L^1$ potentials (Proposition 4.2) we then have $w > 0$ almost everywhere in $\Omega$. Denoting by $v$ the solution of the Dirichlet problem in Proposition 4.1 with datum $(H(w),0)$, the function $w - v$ is nonnegative. By Proposition 3.1 applied to the solution $w-v$ and test function $u_1$, we get

$$0 \leq \int_\Omega \widehat{w-v}\,d(-\Delta u_1 + Vu_1) = \int_\Omega (-H(w))u_1 - \int_{\partial\Omega} \chi_A \frac{\partial u_1}{\partial n}\,d\sigma = -\int_\Omega H(w)u_1,$$

where the last equality follows from the fact that $\partial u_1/\partial n = 0$ on $A$. Thus, the integral in the right-hand side is nonpositive, while the integrand is nonnegative, hence $H(w)u_1 = 0$ almost everywhere in $\Omega$. Since $H(w) > 0$ almost everywhere in $\Omega$, we then have $u_1 = 0$ almost everywhere in $\Omega$. This contradicts the nontriviality of $u_1$. Thus, (6.1) holds.

*Step* 2. We prove the theorem under the additional assumption that the measures $\Delta u_1$ and $\Delta u_2$ are diffuse but not necessarily finite in $\Omega$.

Take a nonzero finite measure $\mu$ in $\Omega$ such that $0 \leq \mu \leq -\Delta u_1 + Vu_1$. In particular, $\mu$ is diffuse and so by Proposition 3.2 the Dirichlet problem (3.2) has a solution $\widetilde{u}_1$ and $\Delta\widetilde{u}_1$ is a finite measure in $\Omega$. Since $V$ is nonnegative, by comparison we have $\widetilde{u}_1 \leq u_1$. By the definition of $\partial u_1/\partial n$ as an essential infimum of normal derivatives over $\mathcal{G}_{u_1}$,

(6.2) $$\frac{\partial u_1}{\partial n} \leq \frac{\partial \widetilde{u}_1}{\partial n} \quad \text{almost everywhere on } \partial\Omega.$$

We next take any $w_2 \in \mathcal{G}_{u_2}$ and apply Proposition 5.2 to this function to get a supersolution $\widetilde{w}_2 \in W_0^{1,2}(\Omega)$ of the Schrödinger operator $-\Delta + V$ such that

(6.3) $$\frac{\partial \widetilde{w}_2}{\partial n} \leq \frac{\partial w_2}{\partial n} \quad \text{almost everywhere on } \partial\Omega.$$

Observe that both $\widetilde{u}_1$ and $\widetilde{w}_2$ satisfy the assumptions of the previous step. Thus,

$$\frac{\partial \widetilde{u}_1}{\partial n} < 0 \quad \text{almost everywhere on } \left\{\frac{\partial \widetilde{w}_2}{\partial n} < 0\right\}.$$

Combining (6.2) and (6.3), we thus get

$$\frac{\partial u_1}{\partial n} < 0 \quad \text{almost everywhere on } \Big\{\frac{\partial w_2}{\partial n} < 0\Big\}.$$

Since this property holds for every $w_2 \in \mathcal{G}_{u_2}$, we obtain (6.1).

*Step* 3. Proof of the theorem completed.

In the general case, it suffices to apply the previous argument to $T_1(u_i)$. Indeed, by Kato's inequality $\Delta T_1(u_i)$ are locally finite measures in $\Omega$. Since $T_1(u_i) \in W_0^{1,2}(\Omega)$, the measures $\Delta T_1(u_i)$ are diffuse with respect to the $W^{1,2}$ capacity. A straightforward variant of Kato's inequality (cf. (5.4) above) implies that $T_1(u_i)$ are supersolutions of $-\Delta + V$. By Step 2, assertion (6.1) above thus applies to $T_1(u_i)$. By (2.2), we have

$$\frac{\partial T_1(u_i)}{\partial n} = \frac{\partial u_i}{\partial n} \quad \text{almost everywhere on } \partial\Omega,$$

and (6.1) for $u_i$ follows. It now suffices to switch the roles between $u_1$ and $u_2$ to conclude. □

*Proof of Theorem 1.* Since the supersolution $u$ is nonnegative, we may replace $V$ by $V^+$, and assume that $V$ is nonnegative. Let $\zeta_0$ be given by Definition 1.1. Since $\partial \zeta_0/\partial n < 0$ almost everywhere on $\partial\Omega$, it suffices to prove that

(6.4) $$\frac{\partial u}{\partial n} < 0 \quad \text{almost everywhere on } \Big\{\frac{\partial \zeta_0}{\partial n} < 0\Big\}.$$

To this end, we apply Proposition 5.2 to any function $w \in \mathcal{G}_{\zeta_0}$ to get a nonnegative supersolution $\widetilde{w} \in W_0^{1,2}(\Omega)$ of the Schrödinger operator $-\Delta + V$ such that $V\widetilde{w} \in L^1(\Omega)$ and

$$\frac{\partial \widetilde{w}}{\partial n} \leq \frac{\partial w}{\partial n} \quad \text{almost everywhere on } \partial\Omega.$$

By Theorem 2, for almost every $x \in \partial\Omega$ we have $\partial u(x)/\partial n < 0$ if and only if $\partial \widetilde{w}(x)/\partial n < 0$. In particular,

$$\frac{\partial u}{\partial n} < 0 \quad \text{almost everywhere on } \Big\{\frac{\partial w}{\partial n} < 0\Big\}.$$

Since this property holds for every $w \in \mathcal{G}_{\zeta_0}$, we have (6.4) and the conclusion follows. □

## 7. EXCEPTIONAL SETS FOR THE HOPF LEMMA

For any negligible compact subset $K \subset \partial\Omega$, we prove that $K$ is the level set $\{\partial u/\partial n = 0\}$ of the normal derivative of a positive smooth solution of

(7.1) $$-\Delta u + Vu = 0 \quad \text{in } \Omega$$

for some $V \in L^1(\Omega; d_{\partial\Omega}\,\mathrm{d}x)$. More generally, given any positive function $\zeta \in C_0^\infty(\overline{\Omega})$, with a normal derivative $\partial \zeta/\partial n$ that possibly vanishes on part of $\partial\Omega$, we find a solution of the equation (7.1), for some $V \in L^1_{\mathrm{loc}}(\Omega)$, whose normal derivative vanishes on a larger subset of $\partial\Omega$ that includes $K$. This is the content of our

**Proposition 7.1.** *Let $\zeta \in C_0^\infty(\overline{\Omega})$ be such that $\zeta > 0$ in $\Omega$. For every compact set $K \subset \partial\Omega$ such that $\mathcal{H}^{N-1}(K) = 0$, there exists $u \in C^\infty(\overline{\Omega} \setminus K) \cap C_0(\overline{\Omega})$ with $0 < u \le \zeta$ in $\Omega$ such that*

  (i) $D^2 u \in L^1(\Omega)$ and $D^2 u / u \in L^1(\Omega; \zeta \, \mathrm{d}x)$,
  (ii) $\partial u / \partial n$ *is well-defined in the classical sense and is continuous on $\partial\Omega$,*
 (iii) $\partial u(x)/\partial n = 0$ *if and only if $x \in K$ or $\partial \zeta(x)/\partial n = 0$.*

*Thus, the function $V := \Delta u / u$ belongs to $L^1(\Omega; \zeta \, \mathrm{d}x)$ by Property (i) above; in particular, $Vu \in L^1(\Omega)$ and*

$$-\Delta u + V u = 0 \quad \text{in } \Omega.$$

It is unclear from our construction whether the upper bound in (1.8) is satisfied by the nonnegative potential $V^+$. We need the following variant of a second-order inequality by Bourdaud [11, Théorème 3]:

**Lemma 7.2.** *Let $H : \mathbb{R} \to \mathbb{R}$ be a convex smooth function such that $H'$ is bounded. If $\zeta \in C_0^1(\overline{\Omega})$ is nonnegative, then for every $\varphi \in C^\infty(\overline{\Omega})$ we have*

$$\|D^2[H(\varphi)]\|_{L^1(\Omega; \zeta \, \mathrm{d}x)} \le C \big(\|D^2 \varphi\|_{L^1(\Omega; \zeta \, \mathrm{d}x)} + \|\nabla \varphi\|_{L^1(\Omega)}\big).$$

*Proof of Lemma 7.2.* In view of the composition formula,

$$D^2 H(\varphi) = H'(\varphi) D^2 \varphi + \nabla[H'(\varphi)] \otimes \nabla \varphi,$$

we only need to estimate the second term in the right-hand side. For every $e \in \mathbb{R}^N$ such that $|e| = 1$, by convexity of $H$ the quantity

$$\big(\nabla[H'(\varphi)] \otimes \nabla \varphi\big)[e, e] = \partial_e[H'(\varphi)] \partial_e \varphi = H''(\varphi)(\partial_e \varphi)^2$$

is nonnegative. Since $\zeta = 0$ on $\partial\Omega$, by integration by parts we get

$$\int_\Omega \big(\nabla[H'(\varphi)] \otimes \nabla \varphi\big)[e, e]\, \zeta = -\int_\Omega H'(\varphi)\big(\partial^2_{e,e}\varphi \, \zeta + \partial_e \varphi \, \partial_e \zeta\big)$$
$$\le \|H'\|_{L^\infty(\mathbb{R})} \big(\|D^2 \varphi\|_{L^1(\Omega; \zeta \, \mathrm{d}x)} + \|\nabla \varphi\|_{L^1(\Omega)} \|\nabla \zeta\|_{L^\infty(\Omega)}\big).$$

This implies the conclusion. □

We also need the following property of the Hausdorff measure $\mathcal{H}^{N-1}$:

**Lemma 7.3.** *Let $K \subset \mathbb{R}^N$ be a compact set. For every $\epsilon > 0$ and every open set $\omega \supset K$, there exists a nonnegative function $\varphi \in C_c^\infty(\omega)$ such that $\varphi > 1$ on $K$ and*

$$\|D^2 \varphi\|_{L^1(\mathbb{R}^N; d_K \, \mathrm{d}x)} + \|\nabla \varphi\|_{L^1(\mathbb{R}^N)} \le C \mathcal{H}^{N-1}(K) + \epsilon,$$

*where $d_K : \mathbb{R}^N \to \mathbb{R}$ denotes the distance to $K$.*

*Proof of Lemma 7.3.* Let $0 < \delta \leq d(K, \partial\omega)/4$, and take finitely many balls $(B_{r_i}(x_i))_{i \in \{1,\ldots,\ell\}}$ that intersect $K$ such that $K \subset \bigcup_{i=1}^{\ell} B_{r_i}(x_i)$,

$$\sum_{i=1}^{\ell} r_i^{N-1} \leq C' \mathcal{H}^{N-1}(K) + \epsilon,$$

and $r_i \leq \delta$ for every $i \in \{1, \ldots, \ell\}$. Given a nonnegative function $\theta \in C_c^\infty(B_2)$ such that $\theta > 1$ on $B_1$, we have the conclusion with

$$\varphi(x) = \sum_{i=1}^{\ell} \theta\Big(\frac{x - x_i}{r_i}\Big).$$

Note that for $x \in B_{2r_i}(x_i)$ we have $d_K(x) \leq 3r_i$. Thus, for every $x \in \mathbb{R}^N$,

$$|D^2\varphi(x)| d_K(x) \leq \sum_{i=1}^{\ell} \frac{3}{r_i} \Big| D^2\theta\Big(\frac{x - x_i}{r_i}\Big)\Big|.$$

A similar pointwise estimate is satisfied by $|\nabla\varphi(x)|$ and we conclude by integration over $\mathbb{R}^N$ and a change of variables in the integral. $\square$

*Proof of Proposition 7.1.* Let $(\epsilon_k)_{k \in \mathbb{N}}$ be a summable sequence of positive numbers. We construct by induction a decreasing sequence of open sets $(\omega_k)_{k \in \mathbb{N}}$ that contain $K$ and a sequence of nonnegative functions $(\varphi_k)_{k \in \mathbb{N}}$ in $C_c^\infty(\omega_k)$ such that $\varphi_k > 1$ on $\omega_{k+1}$ as follows. Take a bounded open subset $\omega_0 \subset \mathbb{R}^N$ that contains $K$ and such that $|\omega_0| \leq \epsilon_0$. Given $\omega_k$, let $\varphi_k \in C_c^\infty(\omega_k)$ be the function given by Lemma 7.3 with open set $\omega_k$ and parameter $\epsilon_k$. We then take an open subset $\omega_{k+1}$ such that

$$K \subset \omega_{k+1} \subset \{\varphi_k > 1\} \quad \text{and} \quad |\omega_{k+1}| \leq \epsilon_{k+1}.$$

We now take a convex smooth function $H : \mathbb{R} \to \mathbb{R}$ such that $H(0) = 1$, $H(t) = 0$ for $t \geq 1$ and $H'$ is bounded. For each $k \in \mathbb{N}$, let

$$\psi_k = H(\varphi_k)\zeta.$$

We have in particular that $\psi_k = 0$ on $\omega_{k+1}$ and $\psi_k = \zeta$ on $\overline{\Omega} \setminus \omega_k$. By the triangle inequality and Lemma 7.2, we have

$$\|D^2\psi_k - D^2\zeta\|_{L^1(\Omega)} \leq \|D^2[H(\varphi_k)]\|_{L^1(\Omega; \zeta\, dx)} + C_1 \|\nabla[H(\varphi_k)]\|_{L^1(\Omega)} + C_2 \|H(\varphi_k) - 1\|_{L^1(\Omega)}$$
$$\leq C_3\big(\|D^2\varphi_k\|_{L^1(\Omega; \zeta\, dx)} + \|\nabla\varphi_k\|_{L^1(\Omega)} + |\omega_k|\big).$$

Since $\zeta \in C_0^\infty(\overline{\Omega})$ and $K \subset \partial\Omega$, we have $\zeta \leq C_4 d_K$ in $\overline{\Omega}$. By the choice of $\omega_k$ and $\varphi_k$ and the assumption $\mathcal{H}^{N-1}(K) = 0$, we deduce that

(7.2) $$\|D^2\psi_k - D^2\zeta\|_{L^1(\Omega)} \leq C_5 \epsilon_k.$$

In particular, the sequence $(D^2\psi_k)_{k \in \mathbb{N}}$ is bounded in $L^1(\Omega)$.

Take

$$u = \sum_{j=0}^{\infty} \frac{1}{2^{j+1}} \psi_j.$$

By construction, we have that $u \in C_0(\overline{\Omega})$, $u$ is smooth in $\overline{\Omega} \setminus K$, and $0 < u \leq \zeta$ in $\Omega$. Moreover, $u$ has a normal derivative given pointwise by

$$\frac{\partial u}{\partial n} = \left(\sum_{j=0}^{\infty} \frac{1}{2^{j+1}} H(\varphi_j)\right) \frac{\partial \zeta}{\partial n}.$$

In particular, $\partial u / \partial n$ is continuous on $\partial \Omega$ and

$$\left\{\frac{\partial u}{\partial n} = 0\right\} = K \cup \left\{\frac{\partial \zeta}{\partial n} = 0\right\}.$$

By the $L^1$ estimate of $D^2\psi_k$, we also have that $D^2 u \in L^1(\Omega)$.

We conclude with the proof that $D^2 u / u \in L^1(\Omega; \zeta \, dx)$. This is based on the pointwise estimate

(7.3) $$\frac{|D^2 u|}{u} \zeta \leq \sum_{j=k}^{\infty} \frac{1}{2^{j-k}} |D^2 \psi_j| \quad \text{on } \omega_k \setminus \omega_{k+1},$$

which follows from the following facts:

(a) $\psi_j = 0$ in $\omega_k$, for every $j < k$,
(b) $u \geq \zeta / 2^{k+1}$ on $\Omega \setminus \omega_{k+1}$, since $\psi_j = \zeta$ on this set, for every $j \geq k+1$.

By (7.3) and the triangle inequality,

$$\frac{|D^2 u|}{u} \zeta \leq \sum_{j=k}^{\infty} \frac{1}{2^{j-k}} |D^2 \psi_j - D^2 \zeta| + 2|D^2 \zeta| \quad \text{on } \omega_k \setminus \omega_{k+1},$$

By this estimate and the fact that $u = \zeta$ in $\Omega \setminus \omega_0$, we have

$$\int_\Omega \frac{|D^2 u|}{u} \zeta = \sum_{k=0}^{\infty} \int_{\omega_k \setminus \omega_{k+1}} \frac{|D^2 u|}{u} \zeta + \int_{\Omega \setminus \omega_0} \frac{|D^2 u|}{u} \zeta$$
$$\leq \sum_{k=0}^{\infty} \sum_{j=k}^{\infty} \frac{1}{2^{j-k}} \int_{\omega_k \setminus \omega_{k+1}} |D^2 \psi_j - D^2 \zeta| + 2 \int_\Omega |D^2 \zeta|.$$

Interchanging the order of summation and using (7.2), we get

$$\sum_{k=0}^{\infty} \sum_{j=k}^{\infty} \frac{1}{2^{j-k}} \int_{\omega_k \setminus \omega_{k+1}} |D^2 \psi_j - D^2 \zeta| = \sum_{j=0}^{\infty} \sum_{k=0}^{j} \frac{1}{2^{j-k}} \int_{\omega_k \setminus \omega_{k+1}} |D^2 \psi_j - D^2 \zeta|$$
$$\leq \sum_{j=0}^{\infty} \int_{\omega_0 \setminus \omega_{j+1}} |D^2 \psi_j - D^2 \zeta| \leq \sum_{j=0}^{\infty} C_5 \epsilon_j < \infty.$$

The proof of the proposition is complete. $\square$

8. POTENTIALS THAT ARE MERELY BOREL FUNCTIONS

Instead of dealing with potentials $V$ in $L^1_{\text{loc}}(\Omega)$, one could wish to work with general Borel functions $V : \Omega \to [0, +\infty]$, but as we explain in this section the counterparts of Theorems 1 and 2 need not be true. The minimization approach that yields a variational solution of the Dirichlet problem

$$\begin{cases} -\Delta\zeta + V\zeta = \mu & \text{in } \Omega, \\ \zeta = 0 & \text{on } \partial\Omega, \end{cases}$$

with datum $\mu \in (W^{1,2}_0(\Omega))'$ can be implemented as in Section 3 above; see [18, 19]. However, since test functions in $C^\infty_c(\Omega)$ need not belong to the minimization class $W^{1,2}_0(\Omega) \cap L^2(\Omega; V \, dx)$, the equation may not be satisfied in the sense of distributions. In this case, the following holds:

**Proposition 8.1.** *Let $V : \Omega \to [0, +\infty]$ be a Borel function and $\mu \in (W^{1,2}_0(\Omega))'$. If $\mu$ is a finite measure in $\Omega$, then the variational solution $\zeta \in W^{1,2}_0(\Omega) \cap L^2(\Omega; V \, dx)$ is such that $V\zeta \in L^1(\Omega)$, $\Delta\zeta$ is a finite measure in $\Omega$ and*

$$\|V\zeta\|_{L^1(\Omega)} + \|\Delta\zeta\|_{\mathcal{M}(\Omega)} \le 3\|\mu\|_{L^1(\Omega)}.$$

*If in addition we have that $\mu \ge 0$ in $\Omega$, then*

$$-\Delta\zeta + V\zeta \le \mu \quad \text{in the sense of distributions in } \Omega.$$

*Proof.* For every $k \in \mathbb{N}$, denote by $\zeta_k$ the minimizer of the functional $E_k$ associated to the bounded potential $V_k := T_k(V)$. In this case, the equation

$$-\Delta\zeta_k + V_k\zeta_k = \mu$$

is satisfied in the sense of distributions and

(8.1) $$\|V_k\zeta_k\|_{L^1(\Omega)} \le \|\mu\|_{L^1(\Omega)}.$$

Thus, $\Delta\zeta_k$ is a finite measure in $\Omega$ and

(8.2) $$\|\Delta\zeta_k\|_{\mathcal{M}(\Omega)} \le \|\mu\|_{\mathcal{M}(\Omega)} + \|V_k\zeta_k\|_{L^1(\Omega)} \le 2\|\mu\|_{\mathcal{M}(\Omega)}.$$

Since $\zeta_k$ is a minimizer of $E_k$ and since $T_k(V) \le V$, for every $k \in \mathbb{N}$ we also have that

$$E_k(\zeta_k) \le E_k(\zeta) \le E(\zeta).$$

We deduce that the sequence $(\zeta_k)_{k \in \mathbb{N}}$ is bounded in $W^{1,2}_0(\Omega)$, whence by this inequality it must converge to the minimizer $\zeta$. By Fatou's lemma, as $k$ tends to infinity in the contraction estimate (8.1) we deduce that $V\zeta \in L^1(\Omega)$. By lower semicontinuity of the norm and estimate (8.2), we also have that $\Delta\zeta$ is a finite measure in $\Omega$.

Observe that if $\mu \ge 0$, then $\zeta_k \ge 0$. By the equation satisfied by $\zeta_k$, as $k$ tends to infinity we deduce from Fatou's lemma that

$$\int_\Omega \zeta(-\Delta\varphi + V\varphi) \le \mu[\varphi] = \int_\Omega \varphi \, d\mu,$$

for every nonnegative function $\varphi \in C_c^\infty(\Omega)$. □

Applying Proposition 8.1 to the positive and negative parts of $\mu$, it follows that there exists a finite measure $\lambda$ in $\Omega$ such that

(8.3) $\quad -\Delta\zeta + V\zeta = \mu + \lambda \quad$ in the sense of distributions in $\Omega$.

This measure $\lambda$ possibly depends on $\mu$ and arises due to the singular character of $V$, but it can vanish even for very singular potentials.

*Example* 8.2. Take the potential $V_\alpha : B_1 \to [0, +\infty]$ defined by

$$V_\alpha(x) = \frac{1}{|x_1|^\alpha}$$

with $\alpha \geq 1$, so that $V_\alpha \notin L^1_{\text{loc}}(B_1)$. We have proved in [37, Proposition 9.2] that for every exponent $\alpha \geq 1$ the Dirichlet problem uncouples in the sense that the variational solution satisfies two independent (homogeneous) Dirichlet problems in $B_1^+$ and $B_1^-$, where

$$B_1^+ = \{x \in B_1 : x_1 > 0\} \quad \text{and} \quad B_1^- = \{x \in B_1 : x_1 < 0\}.$$

Solving separately the Dirichlet problems on $B_1^+$ and $B_1^-$ with $\mu \in L^2(B_1)$ and denoting by $\zeta_+$ and $\zeta_-$ these solutions, the function

$$\zeta := \begin{cases} \zeta_+ & \text{in } B_1^+, \\ \zeta_- & \text{in } B_1^-, \end{cases}$$

belongs to $W_0^{1,2}(B_1) \cap L^2(B_1; V_\alpha \, dx)$ and satisfies

$$\int_{B_1} \zeta(-\Delta\varphi + V_\alpha\varphi) = \int_{B_1} \varphi\mu,$$

for every $\varphi \in C_c^\infty(B_1 \setminus (\partial B_1^+ \cap \partial B_1^-))$. When $\alpha \geq 2$, by Proposition 2.7 this identity actually holds for every $\varphi \in C_c^\infty(B_1)$. Hence,

$$-\Delta\zeta + V_\alpha\zeta = \mu \quad \text{in the sense of distributions in } \Omega,$$

and thus the measure $\lambda$ that satisfies (8.3) is identically zero.

The singularity of $V_\alpha$ in the previous example is so strong that Theorems 1 and 2 are simply false. The reason is that the operator $-\Delta + V_\alpha$ with $\alpha \geq 2$ behaves as if the domain $B_1$ were disconnected, with two connected components $B_1^+$ and $B_1^-$. Indeed, the function $\zeta$ defined above with constant datum $\mu \equiv 1$ satisfies $\partial\zeta/\partial n < 0$ on $\partial B_1 \setminus \{x_1 = 0\}$ by a local application of the classical Hopf lemma. However, the function

$$\widetilde{\zeta} = \begin{cases} \zeta_+ & \text{in } B_1^+, \\ 0 & \text{in } B_1^-, \end{cases}$$

is a also a supersolution for $-\Delta + V_\alpha$, but the normal derivative $\partial\widetilde{\zeta}/\partial n$ is negative only on half of the boundary $\partial B_1$.

ACKNOWLEDGEMENTS

The second author (ACP) was supported by the Fonds de la Recherche scientifique–FNRS under research grants J.0026.15 and J.0020.18. He warmly thanks the Dipartimento di Matematica of the "Sapienza" Università di Roma for the invitation. He also acknowledges the hospitality of the Academia Belgica in Rome.

REFERENCES

[1] A. Ancona, *Une propriété d'invariance des ensembles absorbants par perturbation d'un opérateur elliptique*, Comm. Partial Differential Equations **4** (1979), 321–337.

[2] ______, *Negatively curved manifolds, elliptic operators, and the Martin boundary*, Ann. of Math. (2) **125** (1987), 495–536.

[3] ______, *Elliptic operators, conormal derivatives and positive parts of functions*, J. Funct. Anal. **257** (2009), 2124–2158. With Appendix A by H. Brezis.

[4] ______, *Positive solutions of Schrödinger equations and fine regularity of boundary points*, Math. Z. **272** (2012), 405–427. Erratum: Math. Z. **272** (2012), 429.

[5] G. Anzellotti, *Pairings between measures and bounded functions and compensated compactness*, Ann. Mat. Pura Appl. (4) **135** (1983), 293–318.

[6] ______, *On the extremal stress and displacement in Hencky plasticity*, Duke Math. J. **51** (1984), 133–147.

[7] D. Arcoya and L. Boccardo, *Regularizing effect of the interplay between coefficients in some elliptic equations*, J. Funct. Anal. **268** (2015), 1153–1166.

[8] C. Bandle, V. Moroz, and W. Reichel, *'Boundary blowup' type sub-solutions to semilinear elliptic equations with Hardy potential*, J. Lond. Math. Soc. (2) **77** (2008), 503–523.

[9] M. Bertsch and R. Rostamian, *The principle of linearized stability for a class of degenerate diffusion equations*, J. Differential Equations **57** (1985), 373–405.

[10] M. Bertsch, F. Smarrazzo, and A. Tesei, *A note on the strong maximum principle*, J. Differential Equations **259** (2015), 4356–4375.

[11] G. Bourdaud, *Le calcul fonctionnel dans les espaces de Sobolev*, Invent. Math. **104** (1991), 435–446.

[12] H. Brezis and X. Cabré, *Some simple nonlinear PDE's without solutions*, Boll. Unione Mat. Ital. Sez. B Artic. Ric. Mat. (8) **1** (1998), 223–262.

[13] H. Brezis, M. Marcus, and A. C. Ponce, *Nonlinear elliptic equations with measures revisited*, Mathematical aspects of nonlinear dispersive equations (J. Bourgain, C. Kenig, and S. Klainerman, eds.), Ann. of Math. Stud., vol. 163, Princeton Univ. Press, Princeton, NJ, 2007, pp. 55–109.

[14] H. Brezis and A. C. Ponce, *Remarks on the strong maximum principle*, Differential Integral Equations **16** (2003), 1–12.

[15] ______, *Kato's inequality up to the boundary*, Commun. Contemp. Math. **10** (2008), 1217–1241.

[16] H. Brezis and W. A. Strauss, *Semi-linear second-order elliptic equations in $L^1$*, J. Math. Soc. Japan **25** (1973), 565–590.

[17] G. Dal Maso, *On the integral representation of certain local functionals*, Ricerche Mat. **32** (1983), 85–113.

[18] G. Dal Maso and U. Mosco, *Wiener criteria and energy decay for relaxed Dirichlet problems*, Arch. Rational Mech. Anal. **95** (1986), 345–387.

[19] ______, *Wiener's criterion and $\Gamma$-convergence*, Appl. Math. Optim. **15** (1987), 15–63.

[20] G. Dal Maso, F. Murat, L. Orsina, and A. Prignet, *Renormalized solutions of elliptic equations with general measure data*, Ann. Scuola Norm. Sup. Pisa Cl. Sci. (4) **28** (1999), 741–808.

[21] B. Devyver, M. Fraas, and Y. Pinchover, *Optimal Hardy weight for second-order elliptic operator: an answer to a problem of Agmon*, J. Funct. Anal. **266** (2014), 4422–4489.

[22] J. I. Díaz, *On the ambiguous treatment of the Schrödinger equation for the infinite potential well and an alternative via flat solutions: the one-dimensional case*, Interfaces Free Bound. **17** (2015), 333–351.

[23] ______, *On the ambiguous treatment of the Schrödinger equation for the infinite potential well and an alternative via singular potentials: the multi-dimensional case*, SeMA J. **74** (2017), 255–278.

[24] J. I. Díaz, D. Gómez-Castro, J. M. Rakotoson, and R. Temam, *Linear diffusion with singular absorption potential and/or unbounded convective flow: the weighted space approach*, Discrete Contin. Dyn. Syst. **38** (2018), 509–546.

[25] L. Dupaigne, *Stable solutions of elliptic partial differential equations*, Monographs and Surveys in Pure and Applied Mathematics, vol. 143, Chapman & Hall/CRC, Boca Raton, FL, 2011.

[26] L. C. Evans, *Partial differential equations*, 2nd ed., Graduate Studies in Mathematics, vol. 19, American Mathematical Society, Providence, RI, 2010.

[27] L. C. Evans and R. F. Gariepy, *Measure theory and fine properties of functions*, 2nd ed., Textbooks in Mathematics, CRC Press, Boca Raton, FL, 2015.

[28] D. Feyel and A. de la Pradelle, *Topologies fines et compactifications associées à certains espaces de Dirichlet*, Ann. Inst. Fourier (Grenoble) **27** (1977), 121–146.

[29] D. Gilbarg and N. S. Trudinger, *Elliptic partial differential equations of second order*, Grundlehren der Mathematischen Wissenschaften, vol. 224, Springer-Verlag, Berlin, 1998.

[30] M. Grun-Rehomme, *Caractérisation du sous-différential d'intégrandes convexes dans les espaces de Sobolev*, J. Math. Pures Appl. **56** (1977), 149–156.

[31] T. Kato, *Schrödinger operators with singular potentials*, Israel J. Math. **13** (1972), 135–148.

[32] R. Kohn and R. Temam, *Dual spaces of stresses and strains, with applications to Hencky plasticity*, Appl. Math. Optim. **10** (1983), 1–35.

[33] J.-L. Lions and E. Magenes, *Non-homogeneous boundary value problems and applications, Vol. I*, Die Grundlehren der mathematischen Wissenschaften, vol. 181, Springer-Verlag, New York, 1972.

[34] W. Littman, G. Stampacchia, and H. F. Weinberger, *Regular points for elliptic equations with discontinuous coefficients*, Ann. Scuola Norm. Sup. Pisa (3) **17** (1963), 43–77.

[35] M. Marcus and P.-T. Nguyen, *Moderate solutions of semilinear elliptic equations with Hardy potential*, Ann. Inst. H. Poincaré Anal. Non Linéaire **34** (2017), 69–88.

[36] M. Marcus and L. Véron, *Nonlinear second order elliptic equations involving measures*, De Gruyter Series in Nonlinear Analysis and Applications, vol. 21, De Gruyter, Berlin, 2014.

[37] L. Orsina and A. C. Ponce, *Semilinear elliptic equations and systems with diffuse measures*, J. Evol. Equ. **8** (2008), 781–812.

[38] ______, *Strong maximum principle for Schrödinger operators with singular potential*, Ann. Inst. H. Poincaré Anal. Non Linéaire **33** (2016), 477–493.

[39] Y. Pinchover, *Topics in the theory of positive solutions of second-order elliptic and parabolic partial differential equations*, Spectral theory and mathematical physics: a Festschrift in honor of Barry Simon's 60th birthday, Proc. Sympos. Pure Math., vol. 76, Amer. Math. Soc., Providence, RI, 2007, pp. 329–355.

[40] Y. Pinchover and G. Psaradakis, *On positive solutions of the $(p, A)$-Laplacian with potential in Morrey space*, Anal. PDE **9** (2016), 1317–1358.
[41] Y. Pinchover and K. Tintarev, *Existence of minimizers for Schrödinger operators under domain perturbations with application to Hardy's inequality*, Indiana Univ. Math. J. **54** (2005), 1061–1074.
[42] A. C. Ponce, *Elliptic PDEs, measures and capacities. From the Poisson equation to nonlinear Thomas-Fermi problems*, EMS Tracts in Mathematics, Vol. 23, European Mathematical Society (EMS), Zürich, 2016.
[43] A. C. Ponce and N. Wilmet, *Schrödinger operators involving singular potentials and measure data*, J. Differential Equations **263** (2017), 3581–3610.
[44] M. H. Protter and H. F. Weinberger, *Maximum principles in differential equations*, Springer-Verlag, New York, 1984.
[45] G. Stampacchia, *Le problème de Dirichlet pour les équations elliptiques du second ordre à coefficients discontinus*, Ann. Inst. Fourier (Grenoble) **15** (1965), 189–258.
[46] N. S. Trudinger, *On the positivity of weak supersolutions of nonuniformly elliptic equations*, Bull. Austral. Math. Soc. **19** (1978), 321–324.
[47] L. Véron and C. Yarur, *Boundary value problems with measures for elliptic equations with singular potentials*, J. Funct. Anal. **262** (2012), 733–772. With Appendix A by A. Ancona.
[48] Z. X. Zhao, *Green function for Schrödinger operator and conditioned Feynman-Kac gauge*, J. Math. Anal. Appl. **116** (1986), 309–334.

LUIGI ORSINA
"SAPIENZA" UNIVERSITÀ DI ROMA
DIPARTIMENTO DI MATEMATICA
P.LE A. MORO 2
00185 ROMA, ITALY
*E-mail address*: orsina@mat.uniroma1.it

AUGUSTO C. PONCE
UNIVERSITÉ CATHOLIQUE DE LOUVAIN
INSTITUT DE RECHERCHE EN MATHÉMATIQUE ET PHYSIQUE
CHEMIN DU CYCLOTRON 2, L7.01.02
1348 LOUVAIN-LA-NEUVE, BELGIUM
*E-mail address*: Augusto.Ponce@uclouvain.be